\documentclass[a4paper, 11pt]{article}

\usepackage{amsthm}
\usepackage{amsmath}
\usepackage[a4paper]{geometry}

\usepackage{graphicx}
\usepackage[english]{babel}
\usepackage{amsmath}
\usepackage{lineno}
\usepackage{caption}
\usepackage{comment}
\usepackage{subcaption}
\usepackage{amssymb}
\usepackage{mathtools}
\usepackage{mathptmx} 
\usepackage{xcolor}
\usepackage{tikz}
\usepackage{tkz-euclide}
\usetikzlibrary{positioning}
\usepackage[square,numbers,sort]{natbib}
\usepackage{dsfont}
\usepackage{hyperref}

\newcommand{\nacawidth}{0.825\textwidth}
\newcommand{\derive}[2] {{\frac{\partial {#1} }{\partial {#2}}}}
\newcommand{\derd}[2] {{\frac{\vd #1 }{\vd #2}}}

\newcommand{\Z}{\mathcal Z}

\newcommand{\R}[0]{\mathbb{R}}

\newcommand{\abs}[1]{\left | #1\right| }

\newcommand{\norm}[1]{\left \|#1  \right \|}

\newcommand{\set}[2]{ \left \{ #1 ~\middle|~  #2 \right\}}
\newcommand{\sset}[1]{ \{ #1\}}
\newcommand{\rest}[0]{|}
\newcommand{\Leb}{\mathcal{L}}

\newcommand{\vd}{ \mathrm{d}}
\newcommand{\intd}{\,\mathrm{d}}

\newcommand{\e}{\mathrm{e}}
\renewcommand{\epsilon}{\varepsilon}
\renewcommand{\phi}{\varphi}
\newcommand{\skp}[2]{\left \langle {#1}, {#2} \right \rangle}
\newcommand{\bskp}[3][~]{\left [ {#2}, {#3} \right ]_{#1}}

\newcommand{\rskp}[2]{\left(#1\middle)\cdot\middle(#2\right)}

\theoremstyle{definition}
\newtheorem{definition}{Definition}
\newtheorem{conjecture}{Conjecture}
\theoremstyle{theorem}

\theoremstyle{remark}

\theoremstyle{theorem}

\author{Simon-Christian Klein
			\footnote{simon-christian.klein@tu-bs.de}
			}

\date{\today}

\title{Stabilizing DG Methods Using Dafermos' Entropy Rate Criterion: III -- Unstructured Grids}

\begin{document}

	\maketitle
	\abstract{The approach presented in the second installment of this series is extended to multidimensional systems of conservation laws that are approximated via a Discontinuous Galerkin method on unstructured (triangular) grids. Special attention is paid to predicting the entropy dissipation from boundaries. The resulting schemes are free of tunable viscosity parameters and tested on the Euler equations. The trinity of testcases is the spreading of thermal energy from a point source, transsonic and supersonic flows around airfoils, and supersonic air inlets.
}
	\section{Introduction} \label{sec:intro}
	Solvers for hyperbolic systems of conservation laws
\begin{equation} \label{eq:mdcl}
	\derive{u}{t} + \sum_{l=1}^n \derive{f_l(u(x, t))}{x_l} = 0, \quad f_l(u): \R^k \to \R^k
\end{equation}
in $n$ space dimensions and with $k$ conserved variables are one of the backbones of modern computational fluid dynamics. While such solvers are already in use for 40 years the existence and uniqueness theory for the associated equations is largely open, and counterexamples for the uniqueness of weak solutions exist. One can conjecture that these problems can be solved by the introduction of additional balance laws. 
Assume there exists a differentiable convex function $U: \R^k \to \R$ together with flux functions $F_l(u): \R^k \to \R$ satisfying
\[
	\forall l=1, \ldots, n: \quad \derive {U}{u} \derive{f_l}{u} = \derive{F_l}{u} .
\]
 In this case
\[
	\derive{U(u)}{t} + \sum_{l=1}^n \derive{F_l(u(x, t))}{x_l} = 0 
\]
is satisfied for smooth $u$, as shown in \cite{Lax71}. If $u$ is the strong limit of solutions $u_\epsilon$ satisfying
\[
	\derive{u}{t} + \sum_{l=1}^n \derive{f_l(u(x, t))}{x_l} = \epsilon \nabla^2 u, 
\]
it follows in the sense of distributions \cite{Lax71}
\[
	\derive{U(u)}{t} + \sum_{l=1}^n \derive{F_l(u(x, t))}{x_l} \leq 0.
\]
Enforcing such inequlities at the discrete level
\[
	\derd {U(u(x, t))}{t} +  \sum_{l=1}^n D_l F(u) \leq 0
\] enlarges the robustness of schemes considerably \cite{RBDDRKP2021On}. Still there exist counterexamples for the uniqueness of the multidimensional isentropic \cite{CDK2015Global,CK2018NU,CKMS2021NU} and full Euler equations \cite{FKKM2020Oscillatory} even when this additional inequality is enforced. One can therefore ask if other restrictions can lead to uniqueness in the analytical setting and robustness for numerical approximations. Dafermos proposed the entropy rate criterion \cite{Dafermos72} by postulating that the selected weak solution $u$ should dissipate its total entropy
\[
	E_u(t) = \int_\Omega U(u(x, t)) \intd x, \quad \forall t > 0: \derd{E_u(t)}{t} \leq \derd {E_{\tilde u}}{t} 
\]
at least as fast as every other weak solution $\tilde u$. 
Enforcing such a criterion for numerical approximations seems impossible, as turning on artificial diffusion with growing strength allows unbounded dissipation speed combined with unbounded approximation errors. 
Yet, schemes satisfying this additional criterion in a numerical sense were constructed \cite{KS2023EAR,Klein2023StabilizingI,Klein2023StabilizingII}.
Estimating how much entropy dissipation can occur by an exact weak solution in a finite timespan and steering schemes to dissipate this amount of entropy was a key component in those constructions. We will generalize the results from \cite{Klein2023StabilizingII} in this work to multidimensional DG schemes on unstructured triangulations. 
An overview of multidimensional DG schemes is given in section \ref{sec:theory}. Section \ref{sec:impl} explains how Dafermos' entropy rate criterion can be enforced on unstructured triangulations. 
Tests for the compressible Euler equations in high-energy test cases and around two-dimensional aerospace geometries underline the abilities of the method in section \ref{sec:NT}.
	\section{Theory} \label{sec:theory}
	\subsection{Multidimensional DG Schemes}

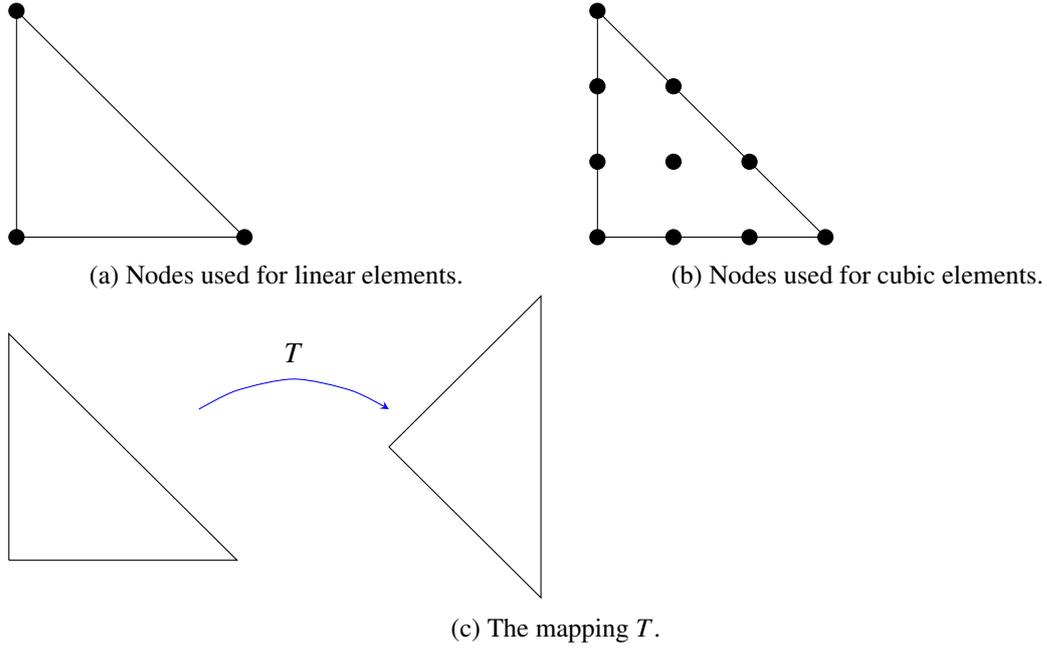
\begin{figure}
		\begin{subfigure}{0.48\columnwidth}
		
		\begin{tikzpicture}[scale=1.0]
			\draw (0.0, 0.0) -- (3.0, 0.0) -- (0.0, 3.0) -- (0.0, 0.0);
			\filldraw (0.0, 0.0) circle (0.1);
			\filldraw (3.0, 0.0) circle (0.1);
			\filldraw (0.0, 3.0) circle (0.1);
		\end{tikzpicture}
		\caption{Nodes used for linear elements.}
	\end{subfigure}
	\begin{subfigure}{0.48\columnwidth}
		\begin{tikzpicture}[scale = 1.0]
			\draw (0.0, 0.0) -- (3.0, 0.0) -- (0.0, 3.0) -- (0.0, 0.0);
			\filldraw (0.0, 0.0) circle (0.1);
			\filldraw (1.0, 0.0) circle (0.1);
			\filldraw (2.0, 0.0) circle (0.1);
			\filldraw (3.0,0.0) circle (0.1);
			\filldraw (2.0, 1.0) circle (0.1);
			\filldraw (1.0, 2.0) circle (0.1);
			\filldraw (0.0, 3.0) circle (0.1);
			\filldraw (0.0, 2.0) circle (0.1);
			\filldraw (0.0, 1.0) circle (0.1);
			\filldraw (1.0, 1.0) circle (0.1);
		\end{tikzpicture}
		\caption{Nodes used for cubic elements.}
	\end{subfigure}
	\begin{subfigure}{0.98\columnwidth}
	\begin{tikzpicture}
		\draw (0.0, 0.0) -- (3.0, 0.0)-- (0.0, 3.0) -- (0.0, 0.0);
		\draw (5.0, 1.5) -- (7.0, 3.5) -- (7.0, -0.5) -- (5.0, 1.5);
		\draw [blue, -stealth] plot [smooth] coordinates {(2.5,2) (3.0, 2.25) (3.75,2.4) (4.5, 2.25) (5,2)};
		\node at (3.75, 2.75) {$T$};
	\end{tikzpicture}
	\caption{The mapping $T$.}
\end{subfigure}
	\caption{Selected nodes for polynomial degrees $1$ and $3$ and the map from the reference element to the physical element.}
	\label{fig:DGnodes}
\end{figure}
Discontinuous Galerkin schemes are high-order generalizations of Finite Volume methods \cite{Cockburn1989DGI,ShuDGReview}.
We assume that the domain $\Omega$ is decomposed into a set of cells $\Z$. Every cell $Z \in \Z$ has a boundary $\partial Z = \cup_{i=1} \partial_i Z$ consisting of several faces $\partial_i Z$ and a well-defined normal $\eta_i$ on every face $\partial_i Z$. In our tests these will be straight-faced triangles without any loss of generality. An approximation space can be build out of polynomials 
\[
	u(x, t)  = \sum_{k=0, l=0}^{k + l \leq p} c_{kl}(t) x_1^k x_2^l
\]
with total degree of less than or equal $p$ on every triangle. We will not assume or enforce any continuity of these ansatz functions between cells, hence the name Discontinuous Galerkin method. Therefore, the ansatz function is defined by $N = \frac{(p+1)(p+2)} 2$ degrees of freedom. 
We select $N$ nodes $x_l$ in every triangle and will from now on deal with the nodal basis $B= \sset{\phi_1, \phi_2, \ldots, \phi_N}$ satisfying
\[
	\phi_l(x_k) = \delta_{kl}.
\]
Our selected nodes are shown in figure \ref{fig:DGnodes}. These correspond to classical first and third order finite elements \cite[example 3.1, 3.3]{Johnson2009FEM} apart from the fact that we place the nodes on the edges on the Gauss-Lobatto points.
Multiplying any component of the conservation law \eqref{eq:mdcl} with a function $v \in V^Z$, integrating over the triangle $Z$ and using the Gauss theorem reveals
\[	
	\begin{aligned}
	&\int_Z v \derd u t + \sum_{l=1}^n \derive{f_l(u(x, t))}{x_l}  \intd V  = 0\\
	 \implies &\int_Z v \derd u t \intd V = \sum_{l=1}^n \int_Z  \derive{v}{x_l}  f_l \intd V - \int_{\partial Z} \rskp{vf(u)}{ \eta} \intd O \\
	 = & \sum_{l=1}^n\int_Z  \derive{v}{x_l}  f_l \intd V - \sum_{m = 1} \int_{\partial_m Z} \rskp{vf(u)}{\eta_m} \intd O.
	 \end{aligned}
\]
We can express this in inner product form
\begin{equation} \label{eq:DGmdbasic}
	\forall v \in V: \quad \skp{v}{\derd u t} = \sum_{l = 1}^n \skp{\derive v {x_l}}{f_l} - \sum_{m=1} \bskp[m]{v}{\rskp{f}{\eta_m}}
\end{equation}
using the inner products
\[
	\skp{u}{v} = \int_Z u v \intd V, \quad \bskp[m]{u}{v} = \int_{\partial_m Z} uv \intd O.
\] 
The weak form \eqref{eq:DGmdbasic} can be converted into a numerical method in two steps. We focus here on nodal Discontinuous Galerkin methods \cite{HW2008DG}, i.e. our approximation of $u$ is saved as a set of nodal values and the corresponding function is defined using the nodal basis $B$. The inner products can be evaluated exactly on the space $V$ using the Gramian matrices
\[
	M_{kl} = \skp{\phi_k}{\phi_l}, \quad S^m_{kl} = \skp{\derive{\phi_k}{x_m}}{\phi_l}, \quad B^m_{kl} = \bskp[m]{\phi_k}{\phi_l}.
\]
Their usage boils down to applying them to the vector of nodal values. The Gauß-Lobatto points on the edges give rise to a surface quadrature, but we will not use any mass lumping in this work, i.e. we will always use these matrices for the discretisation of inner products.
These are only approximations if nodal values of a function that does not lie in $V$ are enterd
\[
	\skp{v}{f} \approx \sum_{k, l} v(x_k) M_{kl} f(u(x_l), t).
\]
Our ansatz functions can be discontinuous between different cells. Solutions of hyperbolic conservation laws harbor discontinuities and one is therefore used to assigning fluxes $f$ to points of discontinuity. Such intercell fluxes are called numerical fluxes in the context of finite volume methods \cite{Toro2009Riemann} and it is natural
to use these fluxes to assign fluxes on the surfaces $\partial Z_m$
\[
	f_m^*(x) = f^{\mathrm{num}}(u_i(x), u_o(x), \eta_m(x)) \approx \rskp {f(u_{\mathrm{R}}(u_i, u_o, 0))}{\eta_m}.  
\]
The exact flux is given by evaluating the solution of the Riemann problem \linebreak[4]$u_R(u_l u_r, x/t)$ at $x = 0$ and inserting it into the continuous flux function\cite{Tadmor1984I,Tadmor1984II}. 
If these are inserted in the surface products follows
\[
 \forall v \in V: \quad	\skp{v}{\derd u t} = \sum_{l=1}^n \skp{\derive{v} {x_l}}{f_l} - \sum_{m=1} \bskp [m] {v}{f_m^*}.
\]
This can be rewritten into matrix-vector form as
\begin{equation} \label{eq:DGmethod}
	\derd u t = M^{-1} \left( \sum_{l = 1}^n S^l f_l - \sum_m B^m f_m^*  \right). 
\end{equation}
We are considering unstructured grids. In general every triangle in the grid will be different. Let us denote the \guillemetleft reference triangle\guillemetright{} as the triangle with vertices at $(0, 0), (0, 1), (1, 0)$. If a triangle in the grid is located at $v_0, v_1, v_2$ then
\[
	T: \R^2 \to \R^2, (r, s) \mapsto T(r, s) = v_0 + r(v_1 - v_0) +s(v_2 - v_0)
\]
maps the reference triangle onto the triangle in the grid. A spatial derivative is therefore given by
\[
	\begin{aligned}
	\derive{u(T(r, s), t)}{x_i} &= \derive{u\left(TT^{-1}(x), t\right)}{x_i} = \derive{u(T(r, s), t)}{(r, s)} \derive{T^{-1}}{x_i}\\
								& = \derive{u}{r} \derive{r}{x_i} + \derive{u}{s}\derive{s}{x_i}.
	\end{aligned}
\]
The boundaries $\gamma_i(t) : [0, 1] \to \partial T_i$ of the reference triangle are mapped to the boundaries of the cell $Z$ as
\[
	\partial Z_i = \set{T (\gamma_i(t))}{t \in [0, 1]} 
\]
holds. The surface normals can be calculated by rotating their tangents by $90$ degrees.
Integrals are transformed as
\[
	\int_{T_{\mathrm{ref}}} u((r, s)) \abs{\det\derive{T}{r, s}} \intd (r, s) =  \int_{Z} u(x) \intd x
\]
for volume and 
\[
\begin{aligned}
	&\int_{\partial Z_i} \rskp{v(x, y)}{\eta_i} \intd O\\
	 = &\int_{0}^{1} \rskp{v(T(\gamma_i(s)))}{\eta_i} \norm{\derive{T(\gamma_i(s))}{s}} \intd s \\
	= &\norm{\partial Z_i} \int_{0}^{1} \rskp{v}{\eta_i}  \intd s \\
	= &\norm{\partial Z_i} \int_{\partial T_i} v(T(r, s)) \eta_i \intd O(r, s)
	\end{aligned}
\]
for surface integrals, where $\norm{\partial Z_i}$ denotes the length of the corresponding edge in the eucledian norm. Therefore the following rules to map the operators $M$, $S^1$, $S^2$, $B^1$, $B^2$, $B^3$ from the reference triangle to an arbitrary triangle in the grid hold,
\[
	\begin{aligned}
	M^Z &= \abs{\det{\derive T {(r, s)}}} M^{T_{\mathrm{ref}}}, \\
	S^{Z, m} &= \abs{\det{\derive T {(r, s)}}} S^{T_{\mathrm{ref}}, r}\derive{r}{x_m} + \abs{\det{\derive T {(r, s)}}} S^{T_{\mathrm{ref}}, s}\derive{s}{x_i},\\
	B^{Z, m} &= B^{T_\mathrm{ref}, m} \frac{\norm{\partial Z_m}}{\norm{\partial T_{\mathrm{ref}, m}}}.
	\end{aligned}
\]

	\section{Implementation of the Entropy Rate Criterion} \label{sec:impl}
	We will now devise how the correction of the per cell entropy given in \cite{Klein2023StabilizingII} can be applied on unstructured two-dimensional grids.
We assume that a cubature $w = (w_1, w_2, \dots, w_N)$ with positive weights, exact at least for our ansatz functions
 exists on every cell $Z$ of the domain subdivision. If such a cubature exists at all is an open question for general node sets, and more nodes than basis functions would be required in general \cite{Glaubitz2023Construction}. Their construction is carried out using the projection onto convex sets algorithm in our case \cite{GPR1970POCS}, similar to the method used for the construction of entire multi-dimensional function space summation by parts operators \cite{GKNO2023FSBP}.
 Let \[
 \begin{aligned}
 C &= \set{w \in \R^N}{w_k \geq 0 } \\
 C_k &= \set{w \in \R^N}{\sum_{i=1}^N w_i \phi_k(x_i) = \int \phi_k(x) \intd x}
 \end{aligned}
 \]
 denote the closed convex set of positive cubature weights and a family of closed convex sets. Every member of the family $C_k$ is a set of cubature weights that is exact for basis function $\phi_k$.
 The maps
 \[
 \begin{aligned}
 	P:& \R^N \to C, w_k \mapsto \max(0, w_k) \\
 	P_k:& \R^N \to C_k, w \mapsto w - \frac{\skp{w}{\left(\phi(x_k)\right)_k} - \int \phi_k \intd x}{\skp{\left(\phi(x_k)\right)_k}{\left(\phi(x_k)\right)_k}} \left(\phi(x_k)\right)_k
 	\end{aligned}
 \]
 are least squares projections onto $C$ and $C_k$. If the intersection of the convex sets $C_k$ and $C$ is nonempty converges the iteration
 \[
 	w^{n+1} = P_N \circ P_{N-1} \circ P_{N-2} \circ \ldots \circ P_1 \circ P w^n
 \]
 to a point $w^\infty$ in their intersection \cite{Neumann1950Functional}. It is therefore a positive and exact cubature.
 Using this cubature we can define a discrete per cell entropy.
 \begin{definition} [Discrete per cell entropy]
 	Let the discrete entropy in cell $Z$ be denoted as $E^Z_u$. We define the discrete entropy by
 	\[
		E^Z_u = \sum_{l=1}^N w_l U(u(x_l, t)).
 	\]	
 \end{definition}
 Our next steps will be the development of entropy inequality predictors for several space dimensions, i.e. functionals $\sigma$ that predict the rate of entropy decay possible in a domain $\theta$, and the development of an entropy dissipation mechanism for the discrete per cell entropy defined above.
 Together we will use both to devise a correction direction $\upsilon$ and strength $\lambda$ to enforce that the solution to
 \[
 	\derd u t = M^{-1} \left( \sum_{l=1}^n S^l f_l - \sum_m B^m f_m^*  \right) + \lambda \upsilon
 \]
 satisfies
 \[
 	\sum_{Z\in\Z}\derd{E^Z_u}{t} \leq \sum_{\theta} \sigma^\theta.
 \]

\subsection{Entropy Inequality Predictors in One Dimension}
\begin{figure}
	\begin{subfigure}{0.48\columnwidth}
		\begin{tikzpicture}
			\draw [line width = 1pt] (0.0, 3.0) -- (4.0, 3.0);
			\node at (2.0, 3.5) {\large $u_l$};
			\node at (2.0, 2.5) {\large $u_r$};
			
			\draw [-stealth, line width = 1pt] (1.0, 2.5) -- (1.0, 1.75) node [left] {$ \Delta t$} --(1.0, 1.0);

			\draw [dashed, line width = 1pt] (0.0, 0.0) -- (4.0, 0.0);
			\draw [blue, line width = 1pt] (0.0, 0.75) -- (4.0, 0.75);
			\draw [blue, line width = 1pt] (0.0, -0.75) -- (4.0, -0.75);
			\node at (2.0, 1.0) {\large $u_l$};
			\node at (2.0, -1.0) {\large $u_r$};
			\node at (2.0, 0.2) {\large $u_{lr}$};
			\draw [-stealth, line width = 1pt] (4.0, 0.0) -- (4.0, 0.25) node [right]{\large $c \Delta t$} -- (4.0, 0.75);
			\draw [-stealth, line width = 1pt] (4.25, 3) -- (4.75, 3) node [right]{\large $x$};
			\draw [-stealth, line width = 1pt] (4.25, 3) -- (4.25, 3.5) node [above]{\large $y$};
		\end{tikzpicture}
		\subcaption{We assume that a weak solution to a Riemann problem with a translation invariant initial condition has a translation invariant solution. Such an initial condition is invariant under the scalings $(x, y) \mapsto (x, by)$ and $(x, y) = (bx, y)$. The solution should therefore satisfy $u(x, y, t) = u(cx, by, bt)$ and $u(x, y, t) = u(bx, y, t)$.}
		\label{fig:TransSymmetry}
	\end{subfigure}
	\begin{subfigure}{0.48\columnwidth}
		\begin{tikzpicture}
			\draw [line width = 1pt] (0.0, 0.0) -- (2.0, 1.5);
			\draw [line width = 1pt] (0.0, 0.0) -- (2.0, -2.0);
			\draw [line width = 1pt] (0.0, 0.0) -- (-2.0, -0.5);
			\draw [line width = 1pt] (0.0, 0.0) -- (-1.0, -2.0);
			\draw [line width = 1pt] (0.0, 0.0) -- (-2.0, 1.5);
			
			\draw [red, line width = 1pt, dashed] (0.0, 0.0) circle (1.5);
			\draw [-stealth, line width = 1pt] (0.0, 0.0) --(0.0, 1.0) node [right] {\large $c \Delta t$} -- (0.0, 1.5);
			\node at (0.1, 0.0) [right]{\large $0$};
		\end{tikzpicture}
		\subcaption{Assume cones of constant states intersect at $0$. This solution is invariant under scaling $(x, y) \mapsto (bx, by)$. We assume therefore that the solution is invariant under rescaling $u(x, y, t) = u(bx, by, bt)$. We further assume that, as in the effectively one-dimensional case in figure \ref{fig:TransSymmetry}, the propagation speed can be bounded by a constant $c$. Therefore, the perturbations from the interaction at $0$ will not leave the ball with radius $ct$.}
		\label{fig:ScalSymmetry}
	\end{subfigure}
	\begin{subfigure}{0.98\columnwidth}
		\centering
		\begin{tikzpicture}
			\draw [line width = 1pt] (-4.0, 0.0) -- (4.0, 0.0);
			\draw [line width = 1pt] (-3.5, -1.33) -- (0.5, 4.0);
			\draw [line width = 1pt] (-0.5, 4.0) -- (3.5, -1.33);

			\draw [red, line width = 1pt, dashed] (-2.5, 0.0) circle (1.0);
			\draw [red, line width = 1pt, dashed] (2.5, 0.0) circle (1.0);
			\draw [blue, line width = 1pt] (-1.6, 0.4) -- (1.6, 0.4);
			\draw [blue, line width = 1pt] (-1.6, -0.4) -- (1.6, -0.4);
			
			\draw [stealth-stealth, line width = 1pt] (2.5, 0.0) -- (3.207, 0.707);
			\node at (2.7, 0.5)  {\large $r$};
			\draw [stealth-stealth, line width = 1pt] (0.0, 0.0) -- (0.0, 0.4);
			\node at (0.3, 0.2) {\large $h$};
		\end{tikzpicture}
		\subcaption{Usage of both cases in a triangle. If a triangle, or more general a polygon, in a grid is considered the invariance in \ref{fig:TransSymmetry} and \ref{fig:ScalSymmetry} only holds in local areas where the perturbations from the missing similarity did not propagate at time $t$.}
		\label{fig:LocSymmetry}
	\end{subfigure}
	\caption{Assumptions on weak solutions. One can construct weak solutions that violate these assumptions, but these are not physically relevant solutions from our point of view.}
\end{figure}
We would like to use a technique similar to the HLL inequality predictor presented in \cite{Klein2023StabilizingII} to quantify the maximal possible entropy dissipation. 
There, a fully discrete entropy inequality predictor is a function $\sigma^\theta(u) : \Leb^2 \to \R$  giving a lower bound on the rate of the entropy dissipation in $\theta$. Let us define the dissipation in the spacetime prism $\theta \times [0, t]$ as
\begin{equation}
		u \mapsto \int_\theta \int_0^t \derive U t + \derive F x \intd t \intd x = s^\theta(u, t).
\end{equation}
$\sigma^\theta$ should be a lower bound on the dissipation rate in the sense that
\begin{equation}
	\sigma^\theta(u) \leq \derd{s^\theta(u, t)}{t} \rest_{t = 0}
\end{equation}
holds.
Finding a sharp lower bound for this time-derivative that can be evaluated fast is a hard problem in itself. We therefore restrict our entropy inequality predictors to piecewise smooth initial data. In this case the entropy equality holds in the smooth areas and we only need to consider the jumps.
It was shown that in the one-dimensional case estimates $a_l, a_r$ on the fastest waves to the left and right are enough to give a simple approximation for the entropy dissipation induced by an initial discontinuity. The main idea is to look at the approximate solution given by Harten, Lax and Leer in their seminal paper and calculate its entropy defect \cite{HLL1984HLL,Klein2023StabilizingII}. This leads to the formula
\begin{equation}
	\begin{aligned}
	\sigma^\theta(u_l, u_r) &=(a_r - a_l) U(u_{lr}) + a_l U(u_l) - a_r U(u_r) +F(u_l) - F(u_r), \\
	u_ {lr} &=  \frac{a_r a_l - a_l u_l + f(u_l) - f(u_r)}{a_r - a_l}
	\end{aligned}
\end{equation}
for Riemann initial data.

\subsubsection{Reduction to a One-Dimensional Problem for Piecewise Constant Initial Data}
Classical Finite Volume methods are often developed by devising a one-\linebreak[3]dimensional inter-cell flux and integrating this flux along the edges of the selected finite volume. This neglects the multidimensional nature of the problem, as there are naturally points where the edges of the cell intersect. There, more complex interactions take place and more than two different states will interact \cite{GCLM2022Entropy}. Still, most schemes use only fluxes calculated on every face of a finite volume.
This reduction to a one-dimensional problem is also our goal as we can then lean on the previous theory and we will give some arguments why this is possible.
Let us for now look at piecewise constant initial data on a triangulation $\Z$. 
We conjecture that the weak solutions that develop from this initial conditions satisfy two basic principles
\begin{itemize}
	\item Finite wave propagation speed,
	\item Local rotational, scaling or translation invariance, for small times.
\end{itemize}
The first bullet is not a new restriction \cite{Einstein1905On}. The second bullet is not at all clear a priori. Counterexamples of weak solutions that satisfy the classical entropy inequality but break the translation invariance are known \cite{AKKMM2020Nonuniqueness,KKMM2020Shocks}. These paradox solutions are rejected by us as they are physically unsound.
We will now give formal definitions of the invariance expected. Assume a situation as in figure \ref{fig:TransSymmetry}, i.e. two constant states separated by a line. We expect a translation invariance along the line separating the states, i.e. along the $x$ direction for this solution
\[
	u(x + b, y, t) = u(x, y, t), \quad  \forall b \in \R.
\]
Scaling invariance should also hold in the $x$ and $y$ direction, i.e. 
\[
	u(bx, by, bt) = u(x, y, t), \quad  \forall b > 0,
\]
when the point $(0,0)$ lies on the line of discontinuity.
This scaling invariance can be also demanded for the initial condition shown in figure \ref{fig:ScalSymmetry}, there a couple of states are separated by rays stemming from the point $(0, 0)$. 
This invariance should also hold in a local sense for more general initial data, i.e. piecewise constant on triangles as shown in figure \ref{fig:LocSymmetry}.
Assume we move the coordinate system to the right lower angle of the triangle shown. 
The emanating red domain of dependence will not interact with the domain of dependence of other crossing lines of discontinuity, for example the red one around the lower left angle of the triangle. 
We can therefore require also a local scaling invariance for small times for all points located in this red domain of dependence. 
If the coordinate system is instead centered at the middle of the lower side of the triangle also a local scaling invariance can be expected there if the points lie near $(0, 0)$. 
Also, small translations in the $x$ direction should be possible. 
As remarked above these restrictions are not to be expected of all weak solutions, but of all weak solutions we consider physically relevant.
Denote by $\mathcal{E} = \sset{E_1, E_2, \ldots, E_N}$ a set of all edges between two constant states and by $\mathcal{V} = \sset{V_1, V_2, \ldots, V_K}$ the set of all vertices where several constant states interact as depicted in figure \ref{fig:LocSymmetry}.
Around every edge and every vertex exists a time dependent set $\omega^V_k(t)$ for vertices, or $\omega^E_l(t)$ for edges big enough so that the waves from the contacting states do not leave that set for small times $t$. 
These sets can be chosen self-similar with respect to time and space rescaling, as we assumed a rescaling self similarity in our assumptions. 
Average values on these sets have to be constant, as the parametrization $\omega^V(t) = \frac{t}{t_0} \omega^{V}(t_0)$ holds for sets associated with vertices 
\[
	\begin{aligned}
	\frac{1}{\mu(\omega^V(t))} \int_{\omega^V(t)} u\left(x, y, t\right) \intd V &= \frac{t_0^2}{t^2 \mu(\omega^V(t_0))} \int_{\omega^V(t_0)} u\left(\frac{t}{t_0}x, \frac{t}{t_0}y, t\right) \det \begin{pmatrix} \frac{t}{t_0} & 0 \\ 0 & \frac{t}{t_0}\end{pmatrix} \intd V \\
	&= \frac{1}{\mu(\omega^V(t_0))} \int_{\omega^V(t_0)} u(x, y, t_0) \intd V
	\end{aligned}
\]
under the assumption that the coordinate center lies on the vertex.
Under the assumption that an edge is without loss of generality aligned with the x-axis, it analogously holds
\begin{equation} \label{eq:TimeTrafo}
	\begin{aligned}
	\frac{1}{\mu(\omega^E(t))} \int_{\omega^E(t)} u(x, y, t) \intd V &= \frac{t_0}{t \mu(\omega^E(t_0))} \int_{\omega^E(t_0)} u\left(x, \frac{t}{t_0} y, t\right) \det \begin{pmatrix} 1 & 0 \\ 0 & \frac{t}{t_0} \end{pmatrix} \intd V\\
	& = \frac{1}{\mu(\omega^E(t_0))} \int_{\omega^E(t_0)} u(x, y, t_0) \intd V.
	\end{aligned}
\end{equation}
Given an arbitrary subdivision $\omega_k$ of the domain satisfying 
\[
	\Omega \subset \bigcup_k \omega_k, \quad \mu(\omega_k  \cap \omega_l) = \delta_{kl} \mu(\omega_k),  
\] 
it follows from Jensen's inequality \cite{1976LaxCalculus} for a convex function $U$ the lower bound
\[
	\begin{aligned}
	\int_\Omega U(u(x, y, t)) \intd V &= \sum_k \frac{\mu(\omega_k)}{\mu(\omega_k)} \int_{\omega_k} U(u(x, y, t)) \intd V\\
								 &\geq \sum_k \mu(\omega_k) U\left(\frac{1}{\mu(\omega_k)} \int_{\omega_k} u(x, y, t)\intd V\right). 
	\end{aligned}
\]
If $u$ is constant on every $\omega_k$ "=" holds in the equation above. A subdivision for which $u \rest_{\omega_k}$ is constant is the triangle grid $\Z = \sset{Z_1, Z_2, \ldots}$ at $t=0$. A second subdivision for which it does not hold is the subdivision into 
\[
\begin{aligned}
	\omega_k^V, \quad \omega_l^E, \quad \omega^{Z \setminus (E \cup V)} = Z \setminus \left(\bigcup_k \omega^E_k	\cup \bigcup_l \omega^V_l  \right),
\end{aligned}
\]
where the last family of sets is built out of the missing parts of $\Omega$.
If the first subdivision is used for the total entropy of the initial condition and the second subdivision to bound the total entropy at $t$ from below it follows
\[
	\begin{aligned}
	s(\Omega, t)= & E_u(t) - E_u(0) \\
	\geq &\sum_{k} \mu\left(\omega_k^V\right) U\left(u_k^V\right) + \sum_{l} \mu\left(\omega_l^E\right) U\left(u_l^E\right) + \sum_{m} \mu\left(\omega^Z_m\right) U\left(u^Z\right) \\
		 &-\sum_{Z \in \Z} \mu(Z) U\left(u^Z\right)\\
	  = &\sum_{k} \mu \left(\omega_k^V\right) \left(U\left(u_k^V \right)- \sum_Z \alpha_Z^{V_k} U\left(u^Z\right)\right) \\
	   & + \sum_l \mu\left(\omega_l^E\right) \left(U\left(u_l^E\right) - \sum_Z \alpha_Z^{E_l} U\left(u^Z\right) \right)\\
	   = & \sum_{k} \mu \left(\omega_k^V\right) C_1 + \sum_l \mu\left(\omega_l^E\right) C_2.
	\end{aligned}
\]
Here, the averages of $u$ on $\omega_k^E$ and $\omega_l^V$,
\[
\begin{aligned}
	u^E_k = \frac{1}{\mu\left(\omega_k^E\right)} \int_{\omega_k^E} u(x, y, t) \intd V, & \quad u^V_l = \frac{1}{\mu\left(\omega_l^V\right)} \int_{\omega_l^V} u(x, y, t) \intd V,
	\end{aligned}
\]
were used. The coefficients $\alpha_Z^{V_k}$ and $\alpha_Z^{E_l}$ shall be defined as
\[
	\begin{aligned}
	\alpha_Z^{V_k} = \frac{\mu(Z \cap V_k)}{\mu(V_k)}, & \quad \alpha_Z^{E_l} = \frac{\mu(Z \cap E_l)}{\mu(E_l)}.
	\end{aligned}
\]
These do not depend on $t$ for small $t$, satisfying
\[\begin{aligned}
	\alpha_Z^{V_k} \in [0, 1], & \quad \sum_Z \alpha_Z^{V_k} = 1
	\end{aligned}
\] and
\[
	\sum_k \mu(\omega_k^V) \left(\sum_Z \alpha_Z^{V_k}\right) = \sum_{Z} \mu(Z) - \sum_{m} \mu(\omega^Z_m),
\]
as they split the volume of $\omega_k^V$ into parts lying in the respective cells $Z$. 
For $t \to 0$ the two parts of the sum for edges and vertices behave like the volumes of $\omega_k^V$ and $\omega_k^E$ as $C_1$ and $C_2$ are constant. Taking the derivative of $s(\Omega, t)$  and inserting $0$ it follows, as the size of $\omega_k^V$ is purely quadratic,
\[
	\sigma^\Omega = \derd{s(\Omega, t)}{t}(0) = \sum_l \derd{\mu\left(\omega_l^E\right)}{t} \left(U\left(u_l^E\right) - \sum_Z \alpha_Z^{E_l} U\left(u^Z\right) \right).
\]
Therefore, the lower bound on the dissipation speed at $t=0$ is reduced to the bound derived from the edges of the cells. This bound is equivalent to the length of the edge multiplied by the one-dimensional dissipation bound.

\subsubsection{Reduction to a One-Dimensional Problem for Piecewise Smooth Data.}
\begin{figure}
	\begin{subfigure}{0.48\columnwidth}
	\includegraphics[width=\textwidth]{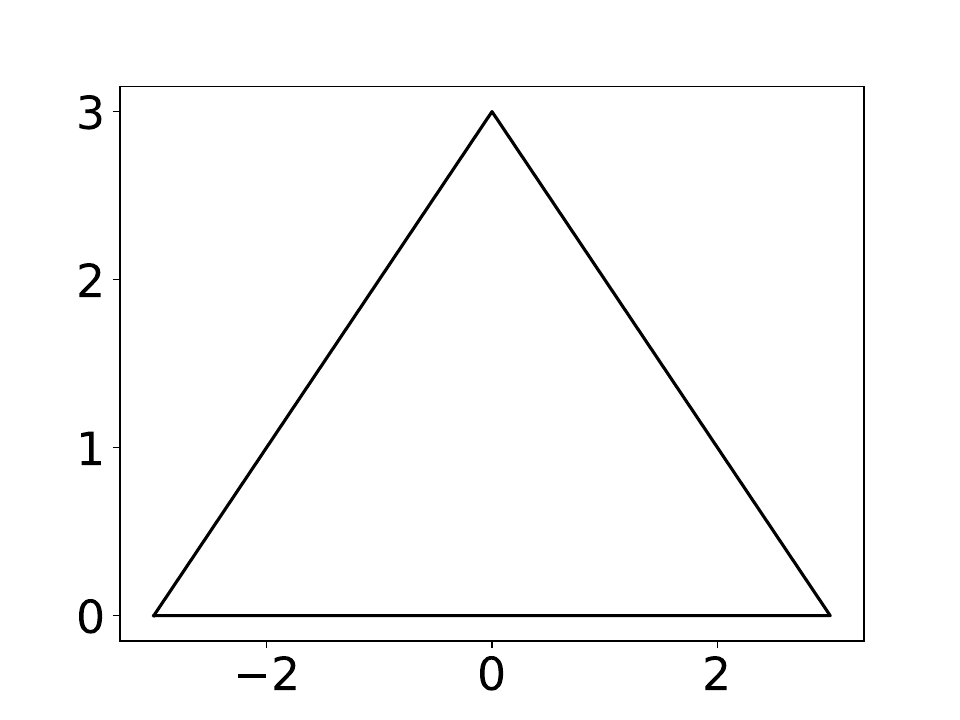}
	\end{subfigure}
\hfill
	\begin{subfigure}{0.48\columnwidth}
		\includegraphics[width=\textwidth]{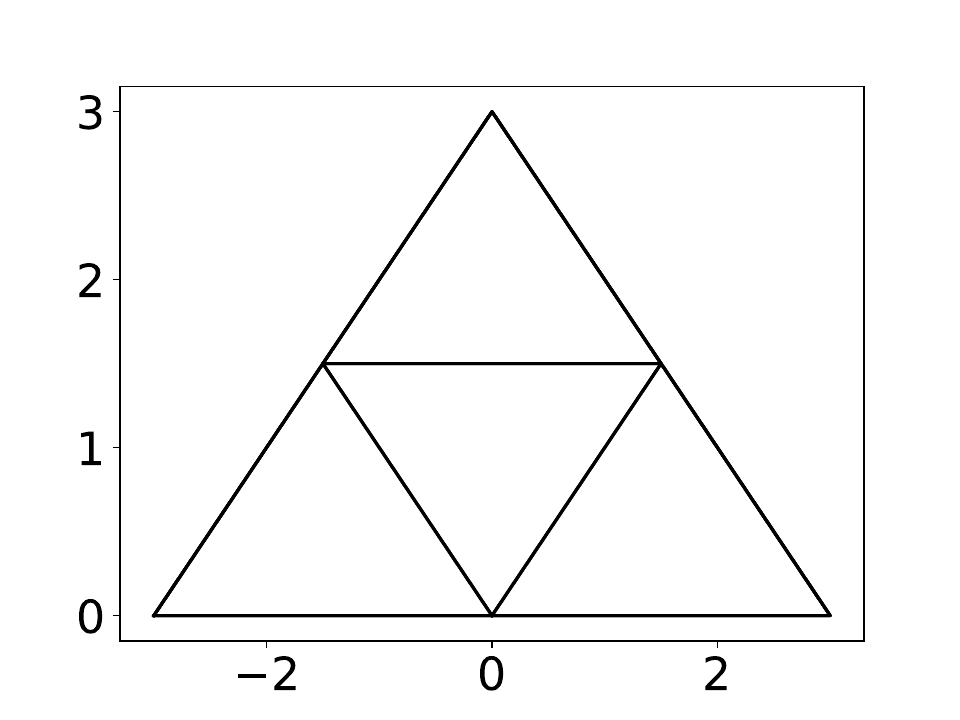}
	\end{subfigure}
	\begin{subfigure}{0.48\columnwidth}
		\includegraphics[width=\textwidth]{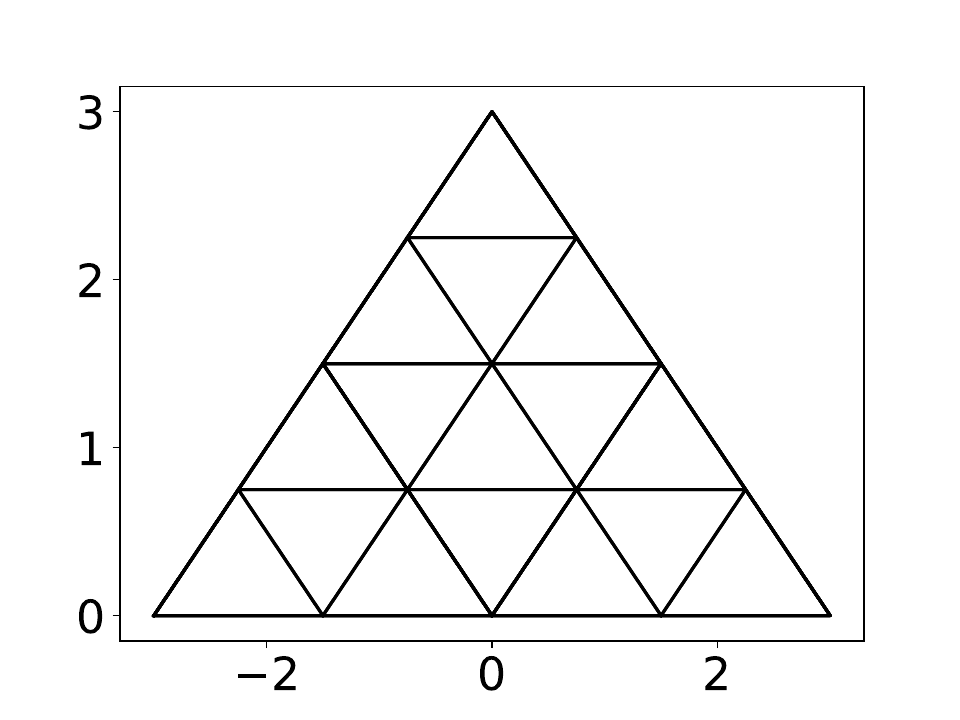}
	\end{subfigure}
\hfill
	\begin{subfigure}{0.48\columnwidth}
		\includegraphics[width=\textwidth]{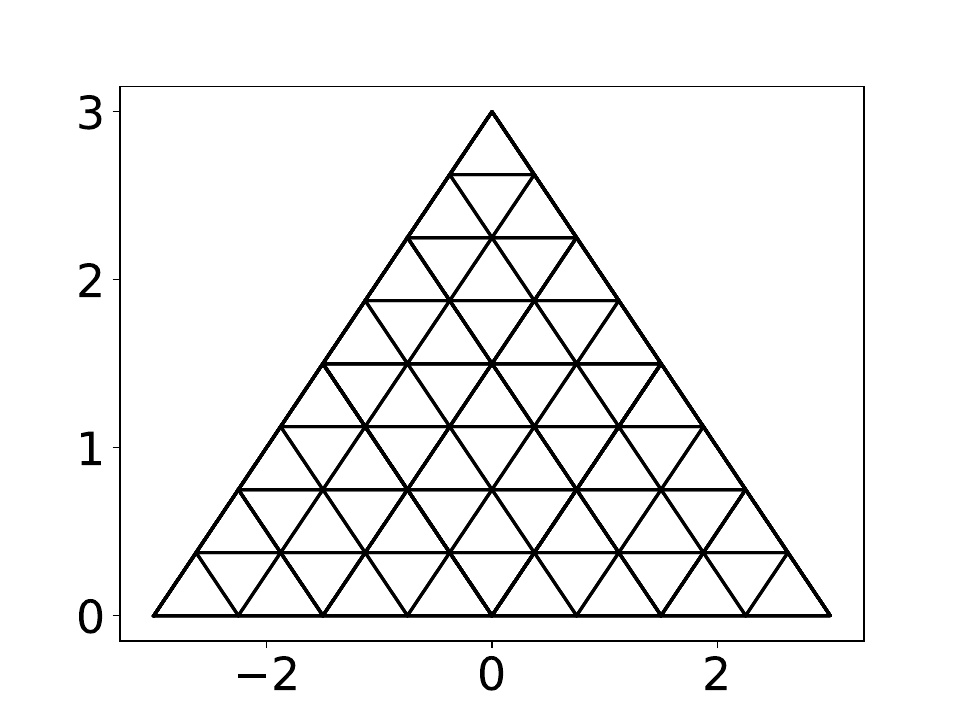}
	\end{subfigure}
	
	\caption{Subdivision of a triangle into finer triangles}
	\label{fig:TrigSubDiv}
\end{figure}

After we explained how one can reduce the problem to one-dimensional entropy inequality predictors for picewise constant data the next step will consist of a generalization to piecewise smooth initial data. Our derivation comes at an additional cost - we will assume that our entropy inequality predictor is a continuous function in the $\Leb^\infty$ norm. 
\begin{conjecture}
Let a sequence of functions $(u_k)_{k=1}^\infty$ converge to a function $u$ in the strong $\Leb^\infty$ sense
\[
	u_k \to u \iff \lim_{k \to \infty} \norm{u_k - u}_{\Leb^\infty} = 0.
\]
Then it follows
\[
	\sigma^\theta(u_k) \to \sigma^\theta(u).
\]
\end{conjecture}
Assume a piecewise smooth solution is given on a set of triangles, i.e. smooth in every triangle but discontinuous across the triangles. 
We can devise a sequence of approximations of $u$ by subdividing every triangle using the sequence depicted in figure \ref{fig:TrigSubDiv}. The following discussion is focused on the reference triangle in figure \ref{fig:TrigSubDiv} for clarity, but mapping this triangle by means of an affine linear map would not interfere with the following arguments.
Every subdivision of this kind divides one triangle into 4 new ones. Therefore the area of every triangle has to be one quarter of the previous one, and every edge is one half as long as before. A triangle that had 6 half edges before is represented by 4 triangles with 9 half edges, as 3 new ones are added to devise the new triangle in the middle of the old one.
Let now $u_{k, l, m}$ be given as piecewise constant approximation of the piecewise smooth function $u$. The function $u_k$ shall be constant on every triangle $\omega_{k, l, m}$ of the $k$-th subdivision of the initial triangle and every of these subdivisions shall be set to the average value of $u$
\[
	u_{k, l, m} = \frac{1}{\mu{\omega_{k, l, m}}} \int_{\omega_{k, l, m}} u(x, y, t) \intd V
\]
 on this subtriangle. Smoothness of $u$ on a subtriangle $\omega_{k, l, m}$, the $m$-th subtriangle of subdivision $k$ of triangle $l$, induces $L$-Lipshitz-continuity. This tells us that the error of the average $u_{k, l, m}$ against $u(x, y)$ is bounded by 
\[
	\begin{aligned}
	\forall (x, y) \in \omega_{k, l, m}: \quad \abs{u_{k, l, m} - u(x, y)} &= \abs{\frac 1 {\omega_{k, l, m}} \int_{\omega_{k, l, m}} u(r, s) - u(x, y) \intd A(r, s)} \\
	&\leq L \underbrace{\max_{(r, s), (x, y) \in \omega_{k, l, m}} \norm{(x, y)^T - (r, s)^T}}_{d_{l, k}},
	\end{aligned}
\]
where the maximum in the equation defines the diameter $d_{l, k}$ of the triangle. This diameter halves with every subdivision, leading to
\begin{equation}
	\label{eq:LinfConv}
	\abs{u_{k, l, m} - u(x, y)} \leq L \frac{d_{l}}{2^k}
\end{equation}
where $d_l$ denotes the diameter of unsubdivided triangle $l$. This shows the (fast) convergence of $u_{k, l, m} \to u$ in the $\Leb^\infty$ norm. Let us split the edges of the $k$-th subdivision into two parts, the edges
\[
	\mathcal{I}_k = \set{E_i \in \mathcal{E}_k}{E_i \text{ is interior}}
\] 
that are interior to a primary triangle, and the edges
\[
	\mathcal{B}_k = \set{E_i \in \mathcal{E}_k}{E_i \text{ is exterior}}
\] 
that are the boundary of a primary triangle. It is easy to see that $\mathcal{B}$ always describes the same total length, while the total length of $\mathcal{I}_k$ grows during every subdivision. We will now count the total length of the interior edges. There are 3 types of interior edges. The ones oriented horizontally, the one oriented in an upwards, and the one oriented in a downwards direction. Focusing on the horizontal ones at subdivision stage $k$ are interior edges lying at the heights 
\[
	h_{k, l} = \frac{l}{2^k}, \quad l = 1, \ldots, 2^k-1
\]
and have the widths
\[
	w_{k, l} = w_0 \frac{2^k - l}{2^k}, \quad l=1, \ldots, 2^k-1.
\]
This amounts to a total length of
\[
	\begin{aligned}
	w_0 \sum_{l = 1}^{2^k-1} \frac{2^k - l}{2^k} = &w_0 \left( 2^k - 1 - \sum_{l=1}^{2^k-1}\frac{l}{2^k} \right) 
	= &w_0 \left(2^k-1 - \frac{\left(2^k -1 \right)2^k}{2^{k+1}} \right)= w_0 \frac{2^k-1}{2}.
	\end{aligned}
\]
We will now sum the total contribution of the interior edges to $\sigma(u_k)$. Please note that in \cite{Klein2023StabilizingII} it was shown that for vanishing $u_l - u_r$ holds $\abs{\sigma(u_l, u_r)} \leq C \norm{u_l - u_r}^2$, and that the jumps at the interior edges vanish as in equation \eqref{eq:LinfConv} because $u(x, y, t)$ is continuous in the interior. One therefore finds
\[
	\abs{\sigma_k(u)\Big\rest_{\mathcal{I}_k}} \leq \sum_{E_i \in  \mathcal{I}_k} \abs{\sigma \left(u_{E_i^-}, u_{E_i^+}\right) \mu(E_i)} \leq 3 \frac{C L^2 d_l^2 w_0 (2^k-1)}{2 (2^k)^2} \overset{k \to \infty}{\to} 0, 
\]
i.e.~the contributions of the interior edges vanish as $u(x, y)$ is smooth in the interior of the master triangle. 
This is also consistent with the entropy equality that holds for smooth $u$ in the analytical context. Looking at the exterior edges one finds 
\[
	\sigma_k(u)\Big\rest_{\mathcal{B}_k} = \sum_{E_i \in \mathcal{B}_k} \sigma\left(u_{E_i^-}, u_{E_i^+}\right) \mu(E_i) \overset{k \to \infty }{\to} \int_{\partial Z_l} \sigma(u_-(x, y), u_+(x, y)) \intd \mathrm{n}(x, y),
\]
i.e.~the corresponding evaluation of the entropy inequality predictor converges to the surface integral on the surface of the master triangle.
We can summarize the previous analysis as follows. One can bound the instantaneous entropy decay for a piecewise smooth initial condition from below under mild assumptions. The bound can be evaluated by integrating a one-dimensional entropy inequality predictor over the boundaries of the triangle. This result can be generalized to other cell structures, like dual cells, by decomposing them into triangles. 

\subsection{First Order FV Schemes Satisfying these Bounds.}
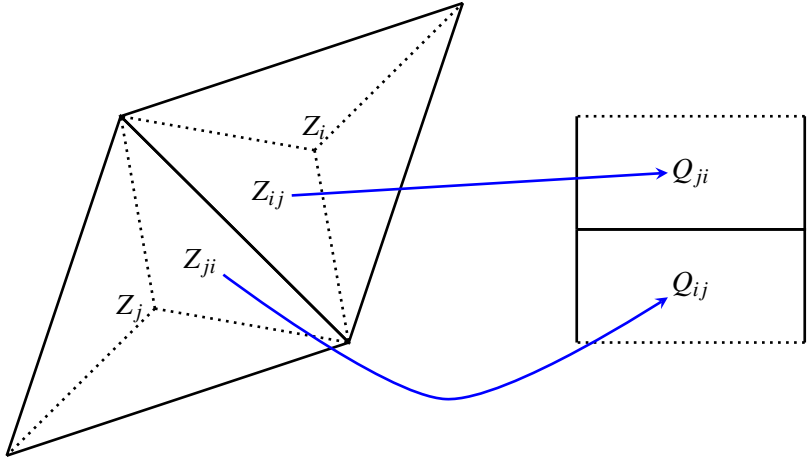
\begin{figure}
	\begin{tikzpicture}[scale = 1.5]
		\draw [line width = 1pt] (0.0, 0.0) -- (3.0, 1.0) -- (1.0, 3.0) -- (0.0, 0.0);
		\draw [dotted, line width = 1pt] (1.3, 1.3) node [left] {$Z_j$} -- (0.0, 0.0);
		\draw [dotted, line width = 1pt] (1.3, 1.3) -- (3.0, 1.0);
		\draw [dotted, line width = 1pt] (1.3, 1.3) -- (1.0, 3.0);
		
		\node at (1.7, 1.7) {$Z_{ji}$};
		
		\draw [line width = 1pt] (4.0, 4.0) -- (3.0, 1.0) -- (1.0, 3.0) -- (4.0, 4.0);
		
		\draw [dotted, line width = 1pt] (2.7, 2.7) node [above] {$Z_i$} -- (3.0, 1.0);
		\draw [dotted, line width = 1pt] (2.7, 2.7) -- (1.0, 3.0);
		\draw [dotted, line width = 1pt] (2.7, 2.7) -- (4.0, 4.0);
		
		\node at (2.3, 2.3) {$Z_{ij}$};
		
		\draw [line width = 1pt] (5.0, 2.0) -- (7.0, 2.0) -- (7.0, 3.0);
		\draw [dotted, line width = 1pt] (7.0, 3.0) -- (5.0, 3.0);
		\draw [line width = 1pt] (5.0, 3.0) -- (5.0, 2.0);
		\node  at (6.0, 2.5) {$Q_{ji}$};
		 \draw [line width = 1pt] (5.0, 2.0) -- (7.0, 2.0) -- (7.0, 1.0);
		 \draw [dotted, line width = 1pt] (7.0, 1.0) -- (5.0, 1.0);
		 \draw [line width = 1pt] (5.0, 1.0) -- (5.0, 2.0);
		 \node at (6.0, 1.5) {${Q_{ij}}$};
		 
		 \draw [blue, -stealth, line width = 1pt] plot [smooth] coordinates {(2.5,2.3) (4.15,2.4)  (5.8,2.5)};
		 \draw [blue, -stealth, line width = 1pt] plot [smooth] coordinates {(1.9,1.6)  (3.85,0.5) (5.8,1.4)};
	\end{tikzpicture}
	\caption{A single step of a scheme formulated on a grid of cells $Z_i \in \Z$ can be rephrased into a set of decoupled schemes for a set of quadrilaterals $Q_{ji}$. After the centroid of a cell $Z_i$ was joined with its vertices results a subdivision of the cell into triangular subcells $Z_{ij}$ that are only adjacent to the subcell $Z_{ji}$ of another cell $Z_j$ and other subcells of $Z_i$. We can assign to every subcell a quadrilateral cell $Q_{ij}$ of equivalent size. The fluxes over the internal, dotted lines, can be set as boundary conditions on the dotted boundaries of the quadrilaterals.}
	\label{fig:TriAv}
\end{figure}
One less obvious question concerns the compatibility of the entropy inequality predictor and the numerical flux function used. Compatibility refers in this case to the question: Does a first order Finite Volume scheme using this respective flux dissipate entropy at least as fast as an entropy inequality predictor demands. This is important, as otherwise the entropy inequality predictor can ask for unattainable rates of entropy descent when the ansatz function in a cell is already constant. The one-dimensional case was already answered positive if a HLL type flux is used in conjunction with an HLL type entropy inequality predictor \cite{Klein2023StabilizingII}.
We can state: An HLL entropy inequality predictor and an HLL flux are compatible on polygonal grids if both use the same speed estimate for $a_l$ and $a_r$. This result is non-obvious as the decoupling of interactions from different boundaries that made the one-dimensional result obvious does not apply in several space dimensions. 
Instead, we will work with a generalization of schemes that can be written in averaging form \cite{Tadmor1984I, Tadmor1984II}. 
In figure \ref{fig:TriAv} two cells $Z_i$ and $Z_j$ out of a polygonal grid are depicted. These can be decomposed into a set of subcells that are only adjacent to one other cell from the original grid. We denote by $Z_{ij}$ the subcell of cell $Z_i$ that is adjacent to cell $Z_j$. 
A state $u_i$ on the cells of the original grid can be spread out to all subcells by interpreting the average value $u_i$ as constant on $T_i$. After one timestep of a first order scheme applied to the subcells $T_{ij}$ their average values can be mapped back to the original grid
\[
	u^{n+1}_i = \sum_{j: T_{ij} \subset T_i} \frac{\mu(T_{ij})}{\mu(T_i)} u^{n+1}_{ij}.
\]
As the scheme on the subcells uses the same fluxes on all edges of the original grid will these averages coincide with averages that were calculated by directly applying the scheme to the cells $T_i$.
Projecting back onto a coarser grid sheds entropy
\[
	U(u_i) \mu(T_i) \leq  \sum_{j: T_{ij} \subset T_i} U(u_{ij}) \mu(T_{ij})
\]
by Jensen's inequality. 
The flux between the subcells of a single macro cell $T_i$ is uniquely determined to be the flux of the average value $f(u_i)$. We can therefore look at an alternative grid of cells $Q_{ij}$ with the same volumes made up of rectangular cells as in figure \ref{fig:TriAv}. Two $Q_{ij}$ and $Q_{ji}$ should be adjacent with an edge of the same length and orientation, and therefore the same flux as between $Z_{ij}$ and $Z_{ji}$. The flux $f(u_i)$ is set as boundary condition on the dotted line in figure \ref{fig:TriAv} and the problem was therefore transformed into one with a translation invariance parallel to the $Q_{ij} \cap Q_{ji}$ edge. Therefore, the one-dimensional argument applies that the entropy inequality prediction is compatible for the cells $Q_{ij}$, and as they have the same surface fluxes and volume, also for $T_{ij}$. As averaging onto the original grid sheds entropy our claim follows.

\subsection{Interplay Between Entropy Inequality Predictors and Boundary Conditions}

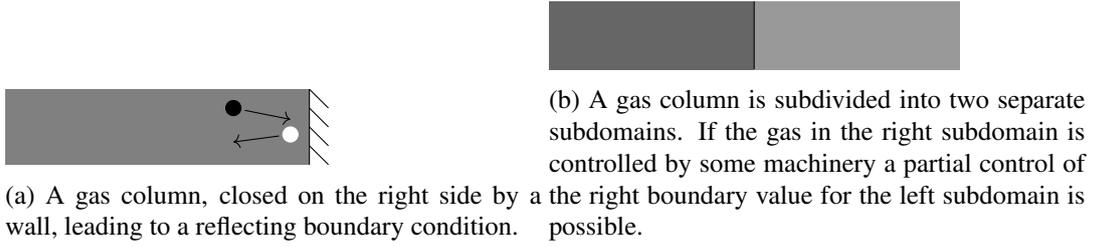
\begin{figure}
	\begin{subfigure}{0.48\columnwidth}
		\begin{tikzpicture}
			\fill [white!50!black] (0.0, 0.0) rectangle (4.0, 1.0);
			\draw (4.0, 0.0) -- (4.0, 1.0);
			\draw (4.0, 0.25) -- (4.25, 0.0);
			\draw (4.0, 0.5) -- (4.25, 0.25);
			\draw (4.0, 0.75) -- (4.25, 0.5);
			\draw (4.0, 1.0) -- (4.25, 0.75);
			\filldraw (3.0, 0.75) circle (0.1);
			\draw [->] (3.15, 0.72) -- (3.75, 0.6);
			\filldraw [white] (3.75, 0.4) circle (0.1);
			\draw [->] (3.6, 0.38) -- (3.0, 0.3);
		\end{tikzpicture}
		\subcaption{A gas column, closed on the right side by a wall, leading to a reflecting boundary condition.}
		\label{fig:bcrefl}
	\end{subfigure}
	\begin{subfigure}{0.48\columnwidth}
		\begin{tikzpicture}[scale=0.9]
			\fill [white!40!black] (0.0, 0.0) rectangle (3.0, 1.0);
			\fill [white!60!black] (3.0, 0.0) rectangle (6.0, 1.0);
			\draw (3.0, 0.0) -- (3.0, 1.0);
		\end{tikzpicture}
		\subcaption{A gas column is subdivided into two separate subdomains. If the gas in the right subdomain is controlled by some machinery a partial control of the right boundary value for the left subdomain is possible.}
		\label{fig:bcdivision}
	\end{subfigure}
	\caption{The two natural boundary conditions for hyperbolic conservation laws that are modelling fluid flows.}
	\label{fig:PBCs}
\end{figure}
\begin{figure}
	\begin{subfigure}{0.48\columnwidth}
		\begin{tikzpicture}
			\filldraw [white!95!black] (0.0, -4.0) -- (2.0, -4.0) -- (2.0, 3.0) -- (0.0, 3.0) -- (0.0, -4.0);
			\draw [->] (0.0, 0.0) node [below] {$0$} -- (0.0, 3.0) node [above] {$t$};
			\draw [->] (-3.0, 0.0) -- (2.0, 0.0) node [right] {$x$};

			\draw (0.0, 0.0) -- (-2.0, 2.0) node [left] {$t a_l = x$};
			\draw (0.0, 0.0) -- (1.0, 2.0) node [right] {$t a_r = x$};
			
			\draw [->] (0.0, -4.0) node [below] {$0$} -- (0.0, -1.0) node [above] {$u$};
			\draw [->] (-3.0, -4.0) -- (2.0, -4.0) node [right] {$x$};
			\draw (-3.0, -2.0) -- (-2.0, -2.0);
			\draw (-2.0, -3.0) -- (1.0, -3.0);
			\draw (1.0, -2.5) -- (2.0, -2.5);
			
		\end{tikzpicture}
		\subcaption{The case $a_l < 0 < a_r$. If the fastest signal speed to the left points into the domain while the fastest signal speed to the right points out of the domain the entropy dissipation will occur in the interior, or exterior, or both areas.}
		\label{fig:BCentrosplit}
	\end{subfigure}
	\hfill
	\begin{subfigure}{0.42\columnwidth}
		\begin{tikzpicture}
			\filldraw [white!95!black] (0.0, -4.0) -- (3.0, -4.0) -- (3.0, 3.0) -- (0.0, 3.0) -- (0.0, -4.0);
			\draw [->] (0.0, 0.0) node [below] {$0$} -- (0.0, 3.0) node [above] {$t$};
			\draw [->] (-1.0, 0.0) -- (3.0, 0.0) node [right] {$x$};
			\draw (0.0, 0.0) -- (1.0, 2.0) node [above] {$t a_l  = x$};
			\draw (0.0, 0.0) -- (2.0, 2.0) node [right] {$t a_r = x$};

			\draw [->] (0.0, -4.0) node [below] {$0$} -- (0.0, -1.0) node [above] {$u$};
			\draw [->] (-1.0, -4.0) -- (3.0, -4.0) node [right] {$x$};
			
			\draw (-1.0, -2.0) -- (1.0, -2.0);
			\draw (1.0, -3.0) -- (2.0, -3.0);
			\draw (2.0, -3.5) -- (3.0, -3.5);
		\end{tikzpicture}
		\subcaption{The case $0 < a_l < a_r$ and vice-versa. If both signal speed estimates lie on one side of the boundary the entropy dissipation can be pinpointed to the interior or exterior of the domain.}
		\label{fig:BCentroout}
	\end{subfigure}
	\caption{Splitting of the HLL approximate solution into internal and external parts. Assume a coordinate system was chosen aligning the normal of the outer boundary with the $x$-axis. The domain lies left of $x = 0$, i.e. the white area, while an outer part is grey. This could be a ghost cell for a reflective boundary condition or the outer cell for a coupling boundary condition.}
	\label{fig:HLLbc}
\end{figure}
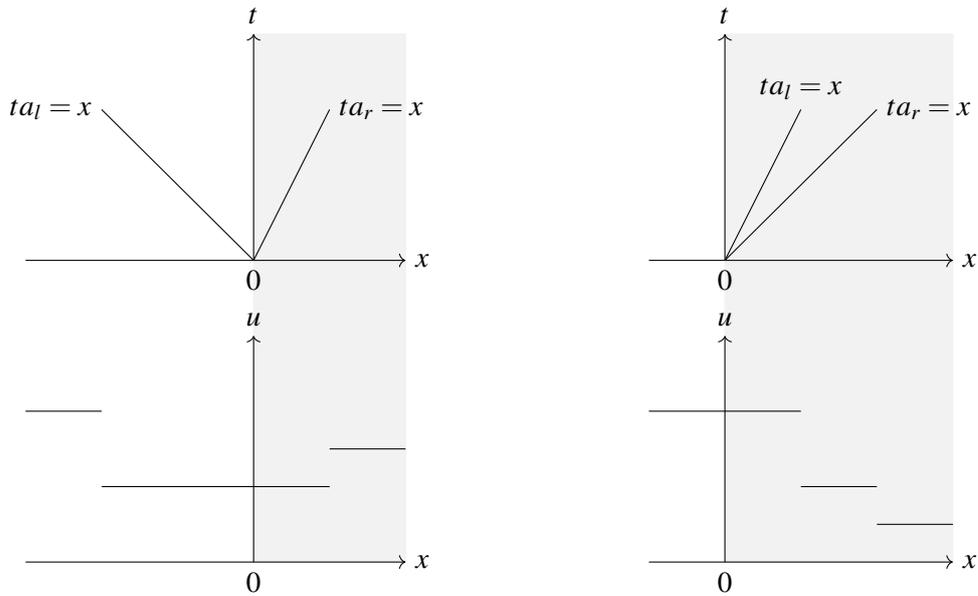

The presentation in the second part of this series focused on periodic domains. While this restriction was only a minor problem for the shock tubes and one-dimensional interactions tested before it is a major concern for what follows. We will therefore describe how reflecting and coupling boundary conditions are handled, especially how a prediction of the residual of the entropy equality can be made at those boundaries. Physical considerations will lead our definitions and this is why only those two boundary conditions, as pictured in figure \ref{fig:PBCs}, will be of interest. Our interpretation follows \cite{CF1948SFSW,Toro2009Riemann,KO2022Entropy}. Far-field boundary conditions do not have a sound physical intuition and will be left out of the discussion. As an intermediate remedy we can only offer to extend simulation domains using growing element sizes up to the point were no waves will reach the far boundaries. This is possible without problems as the method will be able to handle unstructured grids.
The outer ends of such a widely extended domain still correspond to the boundary condition shown in figure \ref{fig:bcdivision}. There the usual entropy inequality predictor can be applied at the boundary values between simulated and non-simulated domains. The interpretation of its value is not as easy as before as only the degrees of freedom in the simulated domain, corresponding to the left one in figure \ref{fig:bcdivision}, can be used to dissipate entropy. If a hypersonic outflow to the right happens in the situation shown should the entropy dissipation happen outside the domain, i.e. the inner degrees of freedom should stay uncorrected as in figure \ref{fig:BCentroout}. This can be taken into account by a slight modification of the HLL entropy inequality predictor. 
Assume after a translation and rotation of the coordinate system is the local normal of the domain boundary parallel to the $x$ axis and $0$ lies on the boundary, the domain lies at negative $x$ values. In other words we changed the coordinate system to coincide with figure \ref{fig:HLLbc}. As explained before, the problem of entropy inequality prediction is a one-dimensional one, and we will from now on deal with the one-dimensional case that has to be integrated over the cell edge afterwards. The entropy inequality prediction on the coupling boundary is now given by
\begin{equation}
	\sigma_{\mathrm{CBC}}(u_l, u_r) = \begin{cases} \sigma(u_l, u_r) & a_r < 0, \\ 
													\sigma(u_l, u_r) & a_l < 0 < a_r, \\ 
														0 & 0 < a_l. \end{cases}
	\end{equation}
	The case in figure \ref{fig:BCentrosplit} is still handled as the case $a_r < 0$, as all dissipation can still happen in the domain, for example if the estimate $a_r$ is just too rough, i.e. too positive.
	Reflective boundary conditions are handled in the same sense, but the external state is always predetermined by the internal state. Density and energy are equal, while the velocity is mirrored \cite{Toro2009Riemann,KO2022Entropy}. This enforces the case $a_l < 0 < a_r$, but the solution should have a mirror symmetry. One can assume that the entropy dissipation is exactly split between external and internal part
	\begin{equation}
		\sigma_{\mathrm{RBC}}(u_l, u_r) = \frac{\sigma(u_l, u_r)}{2}.
		\end{equation}

\subsection{Positive Conservative Filters in Several Space Dimensions}
Dissipation was handled in \cite{Klein2023StabilizingII} using special filters. These were positive, conservative Hilbert-Schmidt operators $\Upsilon$ \cite{LaxFun}.
\begin{definition}[Positive, conservative filter]
	A Hilbert-Schmidt operator 
	\[
	\Upsilon: \Leb^2(\Omega) \to \Leb^2(\Omega), \quad u \mapsto \int_{\Omega} k(x, y) u(y) \intd y
	\] is termed positive if the kernel function $k$ satisfies 
	\[
		\forall x, y \in \Omega \times \Omega: \quad k(x, y) \geq 0.
	\]
	It is called a filter if it averages values of $u$ during its pointwise application
	\[
		\forall x \in \Omega:  \int_\Omega k(x, y) \intd y = 1.
	\]
	The operator is conservative if it satisfies 
	\[
		\forall y \in \Omega: \int_\Omega k(x, y) \intd x = 1.
	\]
\end{definition}
One can show that a conservative filter conserves integral of a conserved quantity \cite{Klein2023StabilizingII}
\[
	\int_\Omega (K u)(x) \intd x = \int_\Omega u(x) \intd x
\]
and that a positive conservative filter dissipates entropy,
\[
	\int_\Omega U((Ku)(x)) \intd x \leq \int_\Omega U(u(x)) \intd x.
\]
A discrete equivalent can be defined with respect to a set of nodes $x_l$ and a cubature rule $w_l$.
\begin{definition}[Positive, conservative, discrete filter]
	Let $\Upsilon \in \R^{n \times n}$ be a square matrix. We call the corresponding mapping
	positive if $ \Upsilon_{kl} \geq 0$ holds, a filter if 
	\[
		\forall k=1, \ldots, n: \quad \sum_{l=1}^n \Upsilon_{kl} = 1
	\]
	is satisfied and conservative if
	\[
		\forall l=1, \ldots, n: \quad \sum_{k=1}^n \Upsilon_{kl} = w_l.
	\]
	\end{definition}
Such an operator is conservative in the sense of not changing the average value with respect to the quadrature rule and one can proof the entropy dissipativity as in the continuous case. An operator $G$ that generates such a positive conservative discrete filter when integrated in time using a forward Euler step $u^{n+1} = u^n + \Delta t Gu^n$  is termed a filter generator \cite{Klein2023StabilizingII}. If the operator $\Upsilon$ is fixed then
\[
	G = \frac{\Upsilon - I} t
\]
is a suitable generator. Let us define an operator $L$ using the base functions $\phi_k$ as
\begin{equation} \label{eq:HeatWOS}
	L = -M^{-1}Q, \quad Q_{kl} = \sum_m \skp{\derive{\phi_k}{x_m}}{\derive{\phi_l}{x_m}}.
\end{equation}
This operator is conservative and generates a filter, but in general not a positive filter. It is related but not equivalent to the heat equation, as integration by parts yields
\begin{equation} \label{eq:HeatWS}
	\skp{u}{\nabla^2 v} = \bskp{u}{\nabla v} - \skp{\nabla u}{\nabla v}.
\end{equation}
The surface terms in equation \eqref{eq:HeatWS} are left out in \eqref{eq:HeatWOS}. We integrate the ODE
\[
	\derd u t = Lu
\]
in time to find the matrix
\[
	C(t) = \e^{tL}.
\]
The operator $C(t)$ is a conservative filter for all times, and
one can show that there exists $t>0$ so that $C(t)$ is positive \cite{Klein2023StabilizingII}. As before, the smallest time $t$ for which this matrix is positive was found via bisection, and used to define the generator
\[
G = \frac{C(t) - I}{t}.
\]
By convexity it holds for $\lambda \in (0, 1]$
\begin{equation}
	E(u + \lambda \Delta t Gu) - E(u) \leq \lambda \left( E(u + \Delta t Gu) - E(u)\right),	
\end{equation}
and we can use this knowledge to bound the dissipativity of $Gu$ from above,
\begin{equation}
	\label{eq:DispSubgradient}
	\skp{\derd U u}{Gu} \leq \frac{E(u + \Delta t Gu) - E(u)}{\Delta t}	\leq 0.
\end{equation}

\subsection{Computation and Application of the Correction}
We have derived all components needed. It is time to explain their usage in multiple dimensions. First a startup procedure is required following these steps:
\begin{enumerate}
	\item Calculate the mass matrix $M$, the stiffness matrices $S^1, S^2$ and boundary matrices $B^i$ for every cell.
	\item Calculate a positive quadrature on the selected node set (abort when impossible).
	\item Calculate the dissipation generator $G$ using a bisection of $\e^{t L}$ for positivity.
\end{enumerate}
Now the time integration is carried out as follows:
\begin{enumerate}
	\item Calculate the uncorrected time derivative $\derd {\tilde u^Z} t$ using the DG method in equation \eqref{eq:DGmethod}.
	\item Evaluate the entropy inequality predictors in every surface node of every cell.
	\item Calculate entropy dissipation estimate by summing over the surface nodes with the surface quadrature.
	\item Calculate the correction direction $\upsilon$ in every cell using the dissipation generator $G$.
	\item Bound the dissipativity of the correction using equation \eqref{eq:DispSubgradient}.
	\item Compute the correction needed for the classical per cell entropy inequality
		\[
		\begin{aligned}
			\sum_k w_k \rskp{\derd U u}{\derd u t + \lambda \upsilon}(x_k) \leq \sum_m \bskp[m]{1}{ F^*} \\
			\implies \lambda_{\mathrm{ED}} = \max\left(0,  \frac{ \sum_m \bskp[m]{1}{ F^*} - \sum_k w_k \rskp{\derd U u}{\derd u t}}{\sum_k w_k \rskp{\derd U u}{\upsilon}}\right).
			\end{aligned}
		\]
	\item Calculate the correction needed to satisfy the entropy dissipation bound 
		\[
			\lambda^{Z_i}_{\mathrm{ER}} = \sum_{Z_o: Z_i \cap Z_o \neq \emptyset} \frac{\sigma(Z_i, Z_o)}{\skp{\derd U u}{Gu}_{w, Z_i} + \skp{\derd U u}{Gu}_{w, Z_o}}.
		\]  
	\item Define the total correction size $\lambda^{Z} = \lambda^Z_{\mathrm{ED}} + \lambda^Z_{\mathrm{ER}}$.
	\item Add the correction to the time derivative from the DG method
	\[
		\derd {u^Z} t = \derd{\tilde u^Z}{t} + \lambda^Z Gu^Z.
	\]
	
\end{enumerate}

	\section{Numerical Tests} \label{sec:NT}
	To show the robustness and accuracy of the two-dimensional entropy fix for DG methods our scheme will be applied to the Euler equations of gas dynamics in two space dimensions given by \cite{Toro2009Riemann}
\begin{equation}
	\derive{u}{t} + \derive {f}{x} + \derive{g}{y} = 0,  \quad u = \begin{pmatrix}
		\rho \\ \rho v_x \\ \rho v_y \\ E
	\end{pmatrix},
	\end{equation}
	with
	\begin{equation}  f = \begin{pmatrix}
	\rho v_x \\ \rho v_x^2 + p \\ \rho v_x v_y \\ v_x (E + p) 
	\end{pmatrix}, \quad g = \begin{pmatrix} \rho v_y \\ \rho v_x v_y \\ \rho v_x^2 + p \\ v_x(E + p) \end{pmatrix}.
\end{equation}
This system of equations is equipped with an entropy 
\begin{equation}
		S(u) = \log(P \rho^\gamma), \quad U(u) = -\rho S,
	\end{equation}
and one entropy flux per coordinate direction \cite{Harten83b}
\begin{equation}
 F(u) = -v_x \rho S, \quad G(u) = -v_y \rho S
\end{equation}

\begin{table}
	\centering
	\begin{tabular}{c | c | c}
		Solver & DDG1 & DDG3 \\
		\hline
		Polynomial Degree & 1 & 3 \\
		Timeintegration & SSPRK33 & SSPRK33 \\
		CFL Number	& 0.5	& 0.1 \\
		Numerical Flux & HLL & HLL \\
		Entropy Inequality Predictor & HLL-based & HLL-based \\
		
		\end{tabular}
		\caption{Properties of the solvers used.}
	\end{table}
\subsection{Accuracy Test}
	High polynomial degrees in a method should also result in high order accurate simulations and these should be almost mesh independent. We will therefore use the linear transport submodel included in the Euler equations for an accuracy test. The initial condition of our test is
	\begin{equation}
		\begin{aligned}
		\Psi(r) &= \exp\left(1 + \frac {-1}{1-r^2}\right) \\
		\rho(x, y, 0) &= \Psi(3 \norm{(x, y) + (1, 0)})^6 + 1 \\
		v_x &= 1, \quad
		v_y = 0, \quad 
		p = 1
		\end{aligned}. 
	\end{equation}
	The Euler equations are solved up to $t = 1.0$ on the domain $\Omega = [-3/2, 3/2] \times [-1/2, 1/2]$. After this time the density variation should have moved from $(-1, 0)$ to $(0, 0)$ and we take the $\Leb^2$ norm of the error to measure the accuracy achieved.
\begin{table}
	\begin{tabular}{c | c | c | c | c }
		triangles & avg. size $A$	& typ.~len.~ $\sqrt A$ 		& $\Leb^2$ Error (p=1) & EOC(p=1) \\
		\hline
		$1577$	& $1.9 \cdot 10^{-3}$ &	$4.36 \cdot 10^{-2} $	& 	$5.88 \cdot 10^{-4}$		&		\\
		$4793$ 	& $6.26 \cdot 10^{-4}$ & $2.5 \cdot 10^{-2} $ 	& 	$1.90 \cdot 10^{-4}$		& 	$2.03$	\\
		$15856$ & $1.89 \cdot 10^{-4}$ & $1.38 \cdot 10^{-2}$	& 	$5.94 \cdot 10^{-5}$ 		&	$2.09$	\\			
	\end{tabular}
	\caption{Results of the accuracy test for $p=1$.}
	\label{tab:convana1}
\end{table} 

\begin{table}
		\begin{tabular}{c | c | c | c | c }
			triangles & avg. size $A$	& typ.~len.~ $\sqrt A$ 		 & $\Leb^2$ Error (p=3) & EOC (p=3)\\
			\hline
			$1577$	& $1.9 \cdot 10^{-3}$ &	$4.36 \cdot 10^{-2} $		& 	$6.85 \cdot 10^{-5}$	&  $~ $\\
			$4793$ 	& $6.26 \cdot 10^{-4}$ & $2.5 \cdot 10^{-2} $ 		&	$1.34 \cdot 10^{-5}$ &   $2.93$\\
			$15856$ & $1.89 \cdot 10^{-4}$ & $1.38 \cdot 10^{-2}$			&	$ 2.60 \cdot 10^{-6}$ &	  $2.76$ \\			
			\end{tabular}
			\caption{Results of the accuracy test for $p=3$.}
			\label{tab:convana3}
		\end{table} 
		The results are reported in the tables \ref{tab:convana1} and \ref{tab:convana3}. Obviously the order of accuracy of the unmodified DG scheme is not degraded by the correction as the entropy inequality prediction vanishes as fast as in the one-dimensional case \cite{Klein2023StabilizingII}.

\subsection{Sedov Blast Wave}
\begin{figure}
			\includegraphics[width=\textwidth]{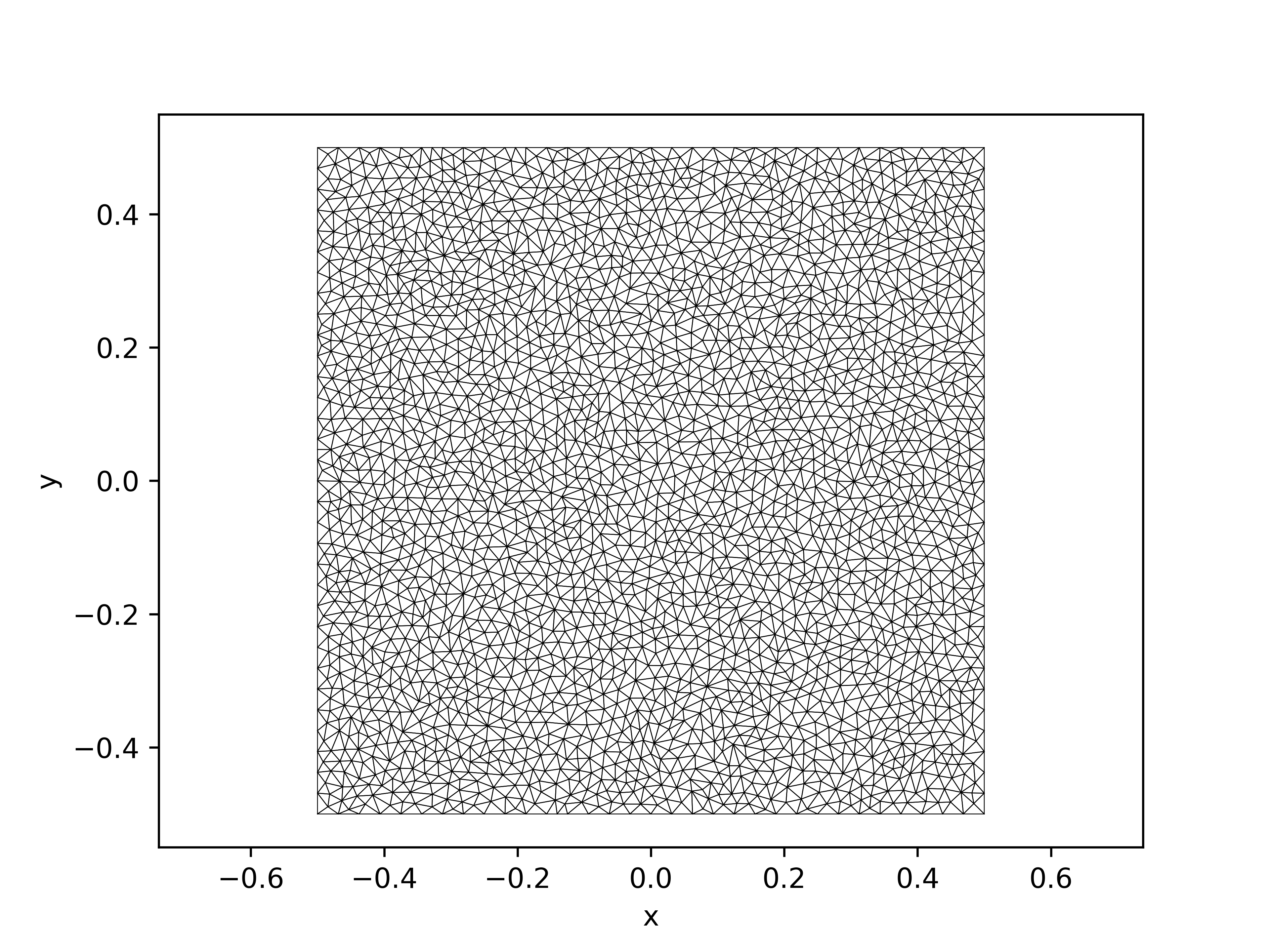}
			\caption{Grid for the Sedov blast wave using 5221 triangles.}
			\label{fig:sedovgrid}
			\end{figure}

\begin{figure}
			\begin{subfigure}{\columnwidth}
				\includegraphics[width=0.9\textwidth]{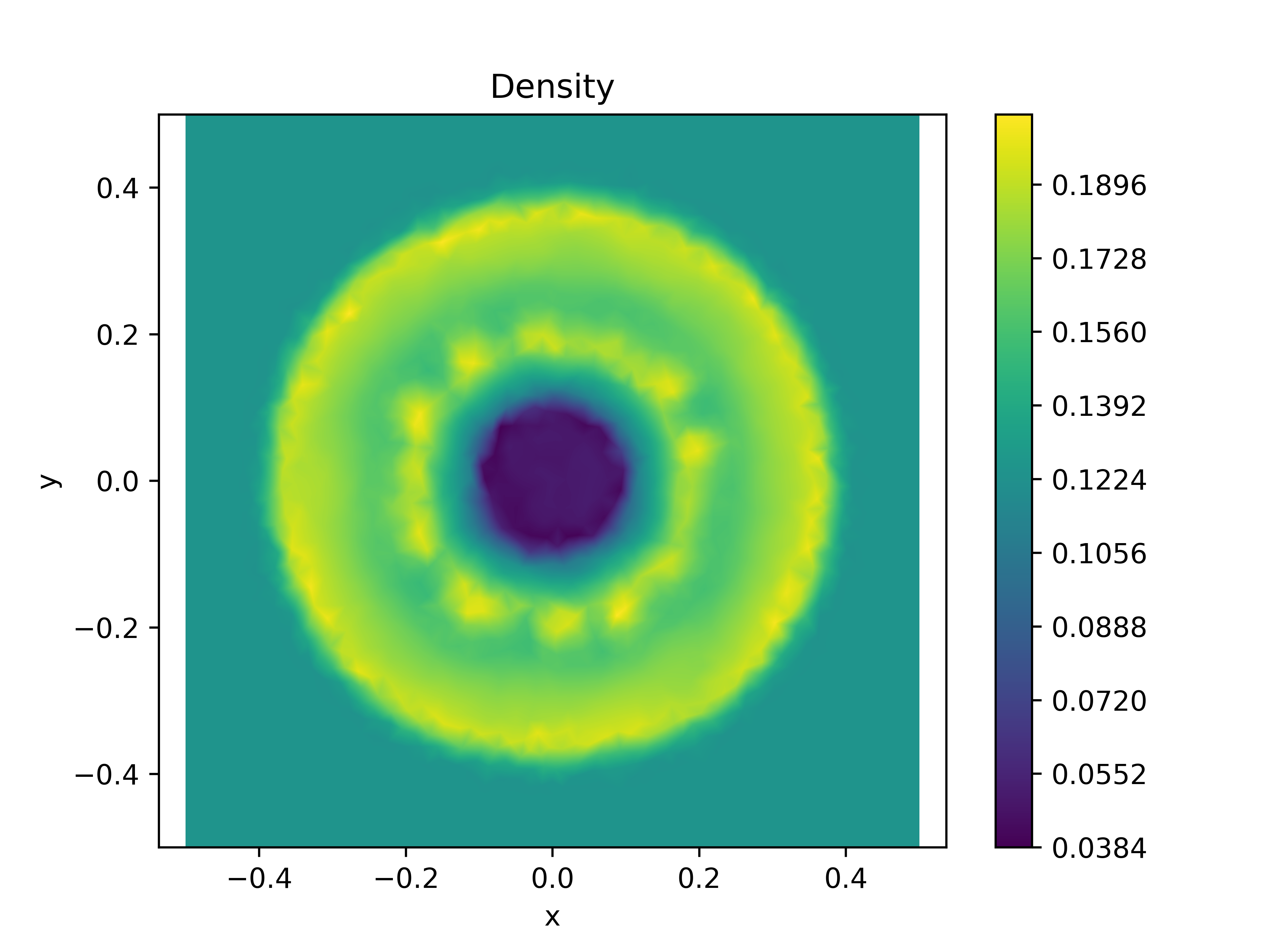}
				\subcaption{Solution at $t=0.2$, $p=1$.}
				\end{subfigure}
				\begin{subfigure}{\columnwidth}
					\includegraphics[width=0.9\textwidth]{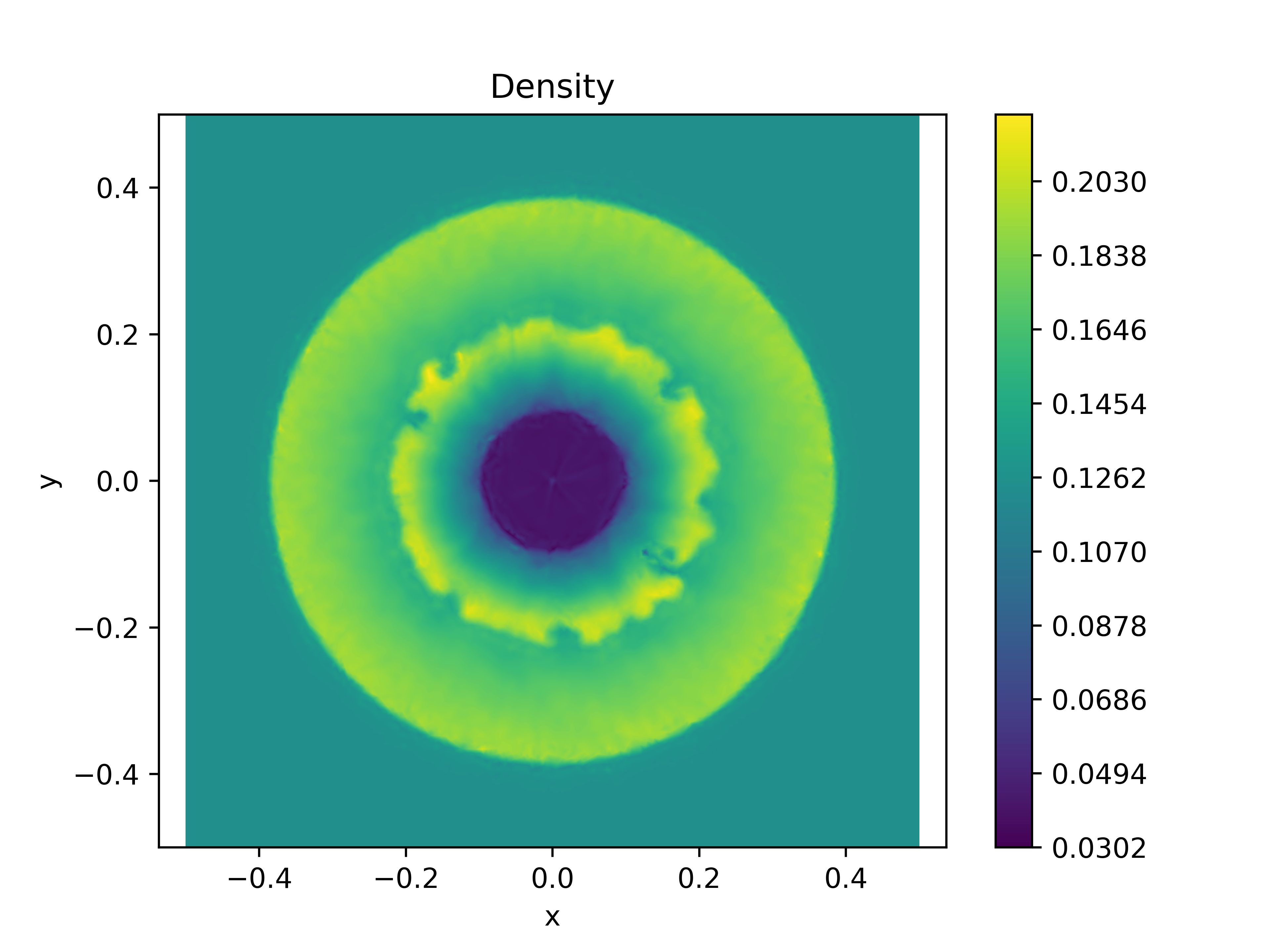}
					\subcaption{Solution at $t = 0.2$, $p=3$.}
					\end{subfigure}
		\caption{Solution to the Sedov blast wave.}
		\label{fig:sedovsol}
	\end{figure}

	The calculation of blast waves can be exemplified by solving Sedov's blast wave problem \cite{KL2007Reduction}. The initial condition
	\begin{equation}
		\begin{aligned}
		u(x, y, 0) &= \begin{cases} u_{\mathrm{in}} & \norm{(x, y)} \leq r \\
					 u_{\mathrm{out}} & \norm{(x, y)} > r  \end{cases}, & \quad& r = 0.08 \\
					 \rho_{\mathrm{in}} &= 1.0, \quad v_\mathrm{in} = 0, &\quad p_{\mathrm{in}}= 1.0\\
					  \rho_\mathrm{out} &= 0.125, \quad v_\mathrm{out} = 0, &\quad p_{\mathrm{out}}= 0.1
					  					 \end{aligned}
	\end{equation}
	has a circular symmetry. The solution is a rotation symmetric blast wave extending from the coordinate origin. Numerical solutions of this problem are often not symmetric. The reason is that the solution consists of two discontinuities traveling outward, and the inner one is a contact discontinuity. This contact is notoriously unstable and is deformed by mesh effects. We use the unstructured triangular mesh in figure \ref{fig:sedovgrid} and the corresponding solutions can be seen in figure \ref{fig:sedovsol}.


\subsection{Forward Facing Step}
\begin{figure}
			\includegraphics[width=\textwidth]{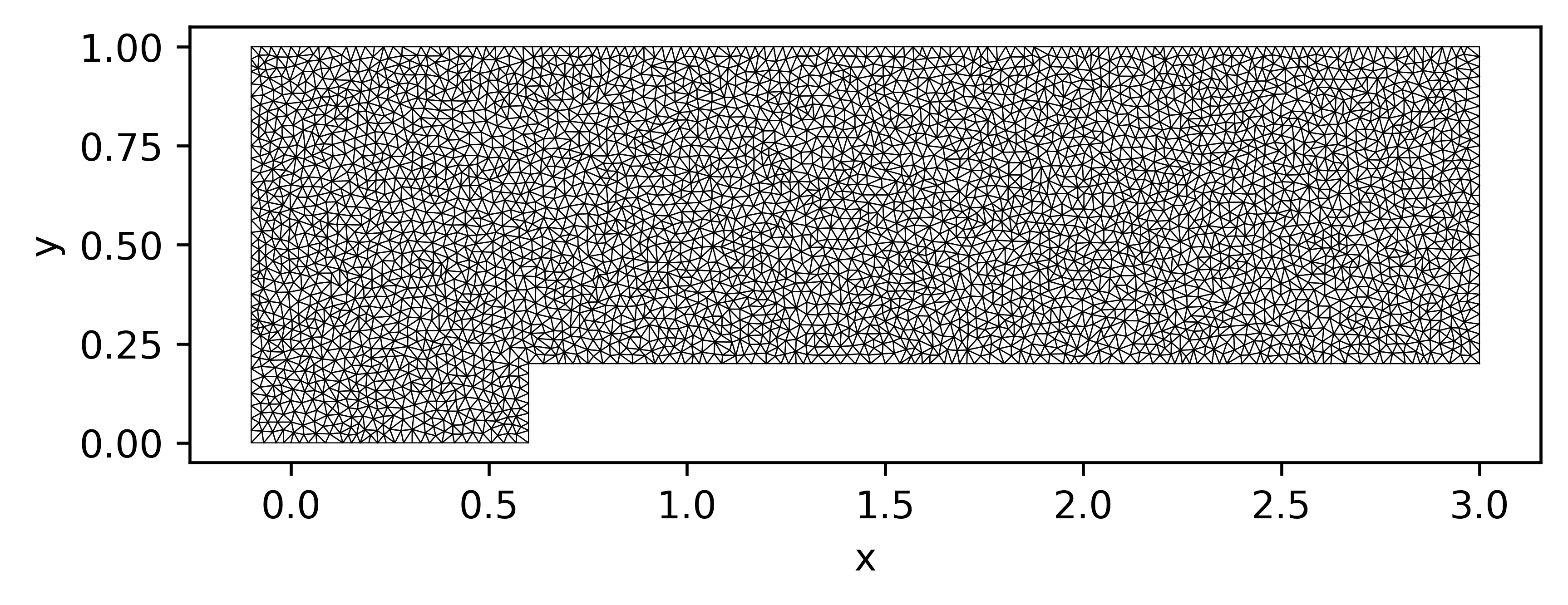}
			\caption{Grid used for the forward facing step problem.}
		\label{fig:ffs}
	\end{figure}

\begin{figure}
			\begin{subfigure}{\columnwidth}
			\includegraphics[width=\textwidth]{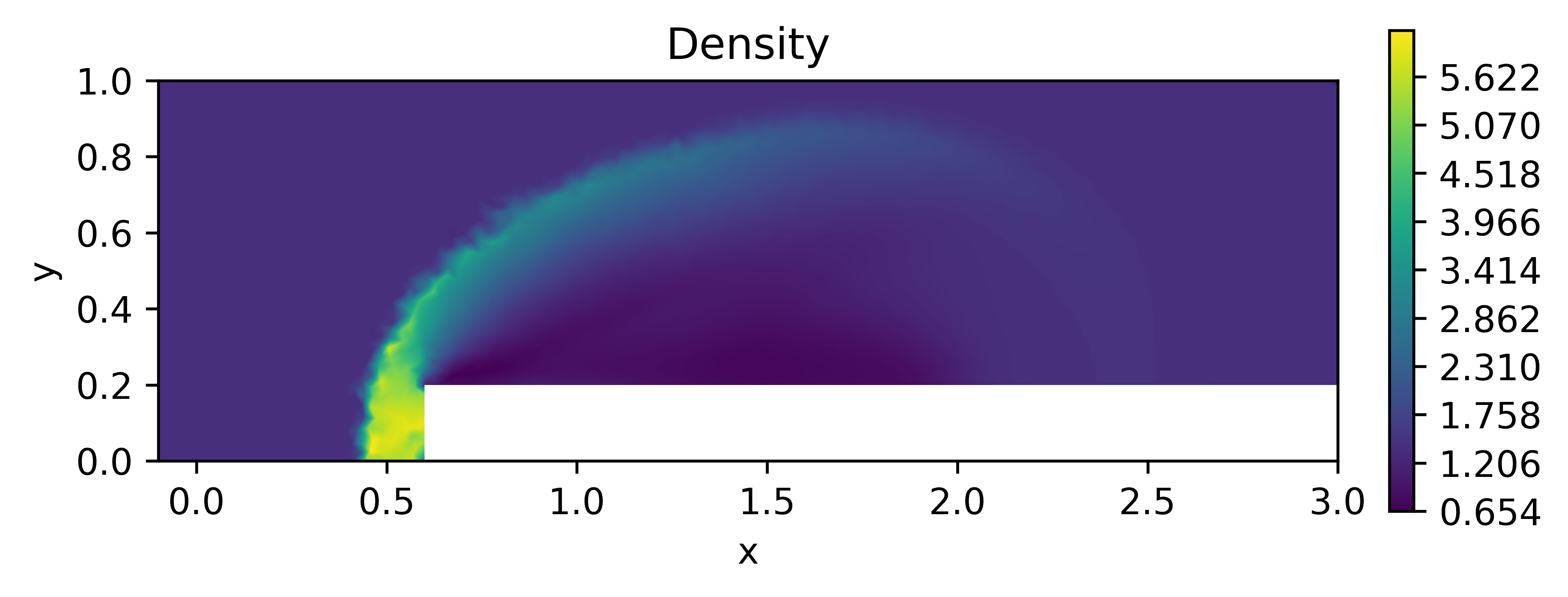}
			\caption{Density at $t=0.5$}
		\end{subfigure}
		\begin{subfigure}{\columnwidth}
			\includegraphics[width=\textwidth]{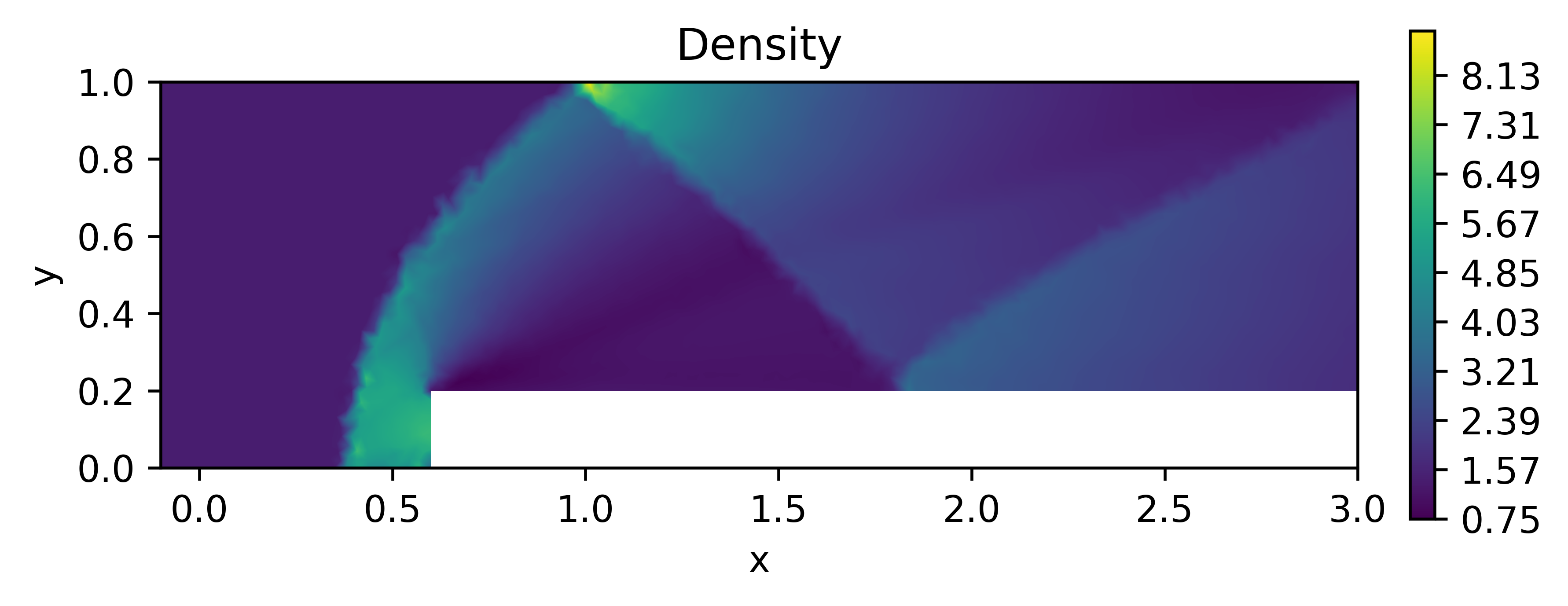}
			\caption{Density at $t=3.0$}
		\end{subfigure}
		\caption{Solution to the forward facing step problem using polynomial degree $p=1$.}
		\label{fig:ffs1}
		\end{figure}

\begin{figure}
		\begin{subfigure}{\columnwidth}
			\includegraphics[width=\textwidth]{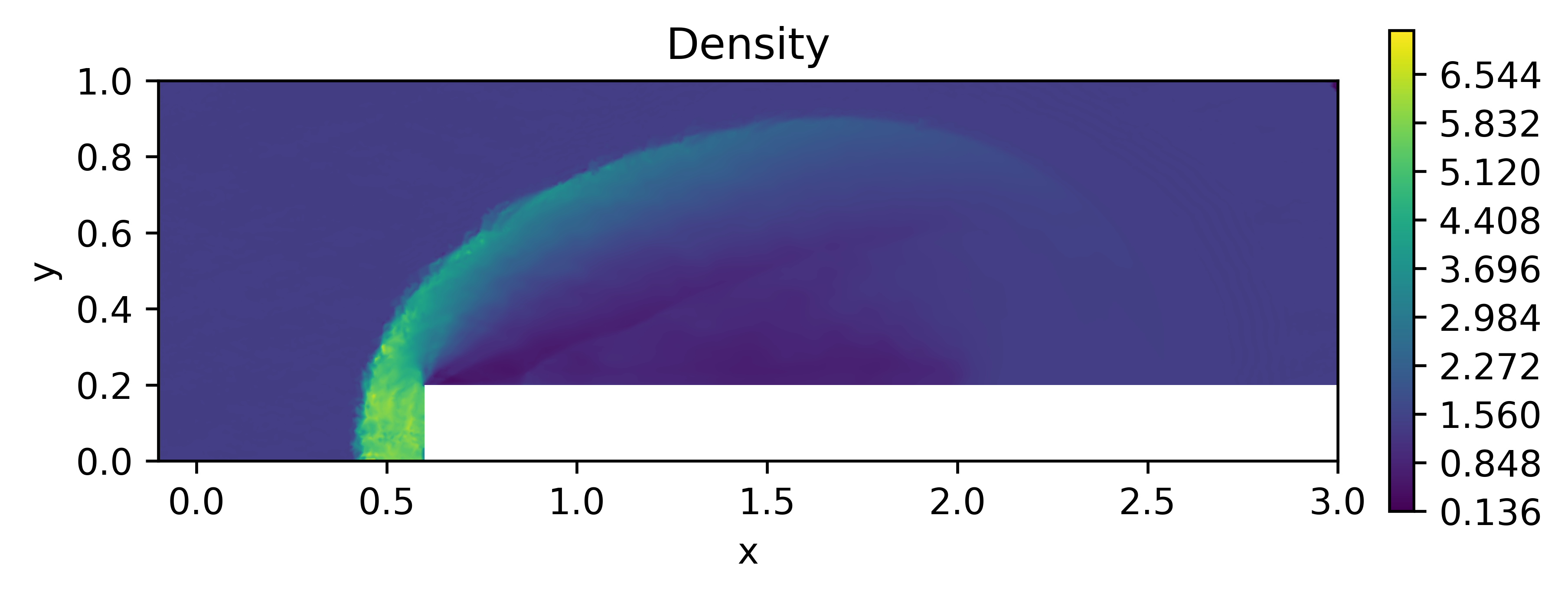}
			\caption{Density at $t=0.5$}
			\end{subfigure}
				\begin{subfigure}{\columnwidth}
				\includegraphics[width=\textwidth]{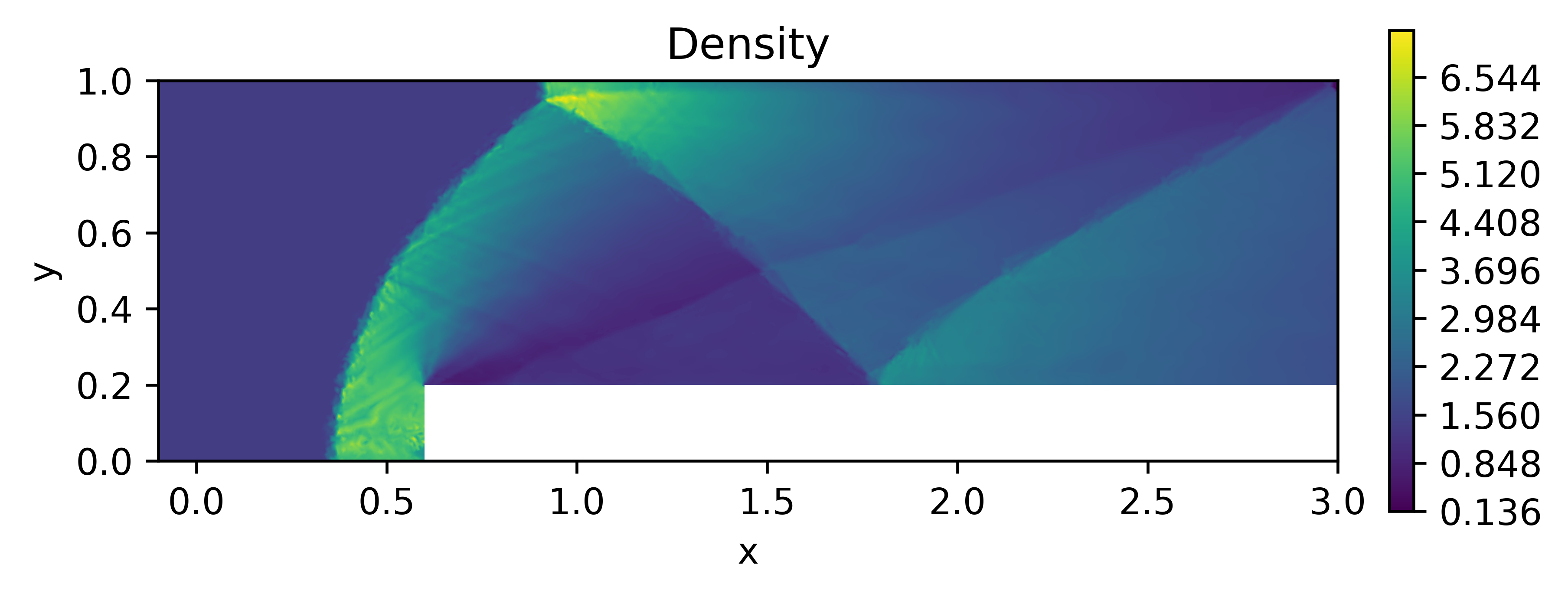}
				\caption{Density at $t=3.0$}
			\end{subfigure}
		\caption{Solution to the forward facing step problem using polynomial degree $p=3$.}
		\label{fig:ffs3}
		\end{figure}

	A model engine inlet can be formed by a step in a channel as for example in figure \ref{fig:ffs}. A channel of length $3$ and height $1$ is equipped with a step at $x = 0.6$ of height $0.2$. This test with a steady mach 3 inflow was proposed in \cite{Emery1968An} and made popular in \cite{WC1984Numerical}. The initial conditions are
	\begin{equation}
		\begin{aligned}
		\rho(x, y, 0) &= \frac{7}{5}, \quad v_x(x, y, 0) = 3 \\
		\quad v_y(x, y, 0) &= 0, \quad p(x, y, 0) = 1.	
		\end{aligned}
	\end{equation}
	This initial condition is used as inflow data on the left end of the tunnel, and the right side should correspond to an outflow. 
	A calculation where degree $p=1$ can be seen in figure \ref{fig:ffs1}, while the same calculation carried out with $p=3$ can be seen in figure \ref{fig:ffs3}. In both cases a grid with the parameters summarized in table \ref{tab:ffsgrid} was used. 
	\begin{table}
		\begin{tabular}{c|c}
			Number of Triangles & 8019 \\
			Number of Nodes		& 8019 \\
			Minimum angle		& 28  \\
			\end{tabular}
			\caption{Table summarizing the parameters of the grid used for the forward facing step problem.}
			\label{tab:ffsgrid}
		\end{table}
	
\subsection{Naca 0012 Airfoil}
\begin{figure}
	\begin{subfigure}{\columnwidth}
		\includegraphics[width=\nacawidth]{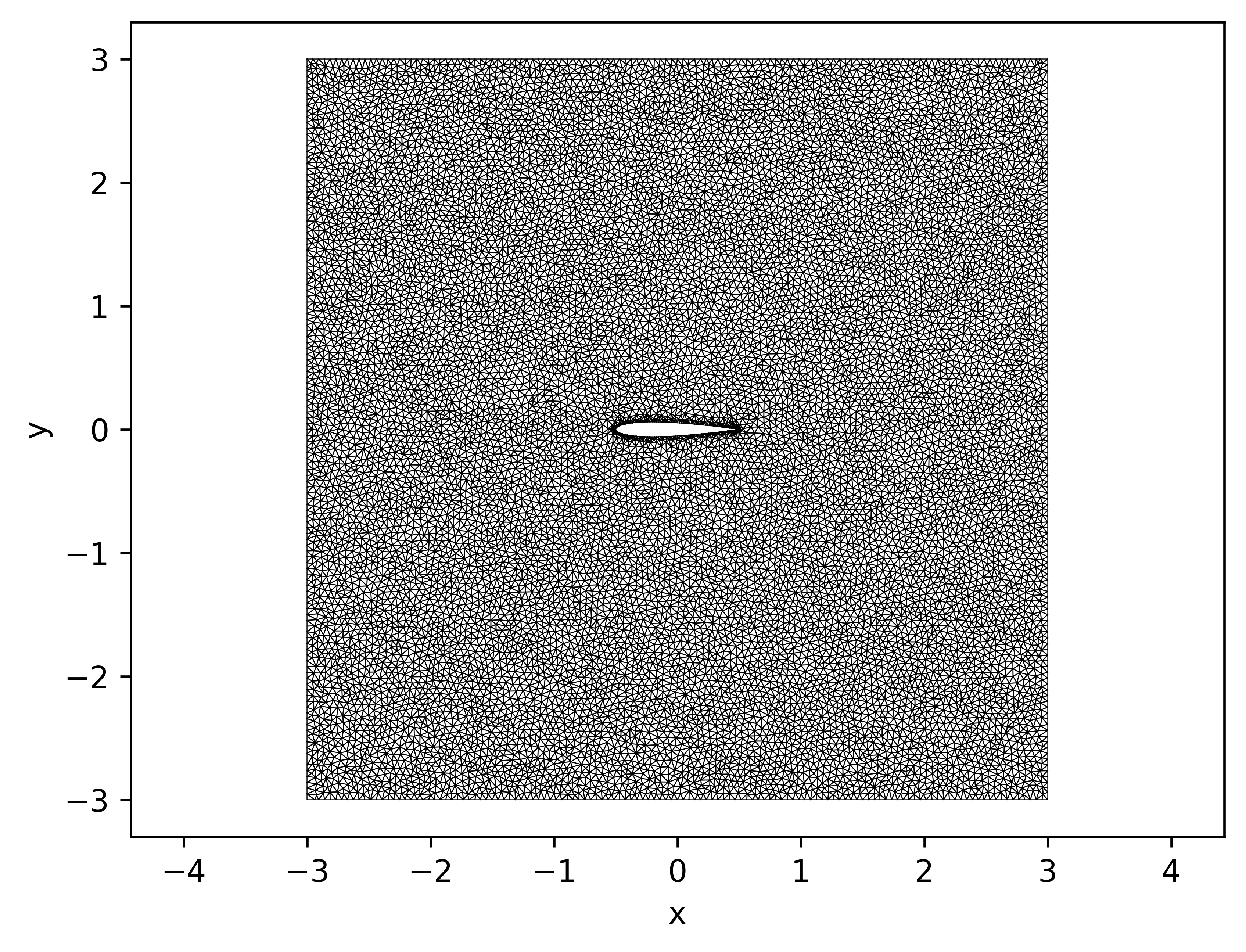}
		\caption{Used grid. In total 27771 triangles were used to triangulate the rectangle $[-3, 3] \times [-3, 3]$ with a maximal area of $0.002$ per triangle}
		\end{subfigure}
	\begin{subfigure}{\columnwidth}
		\includegraphics[width=\nacawidth]{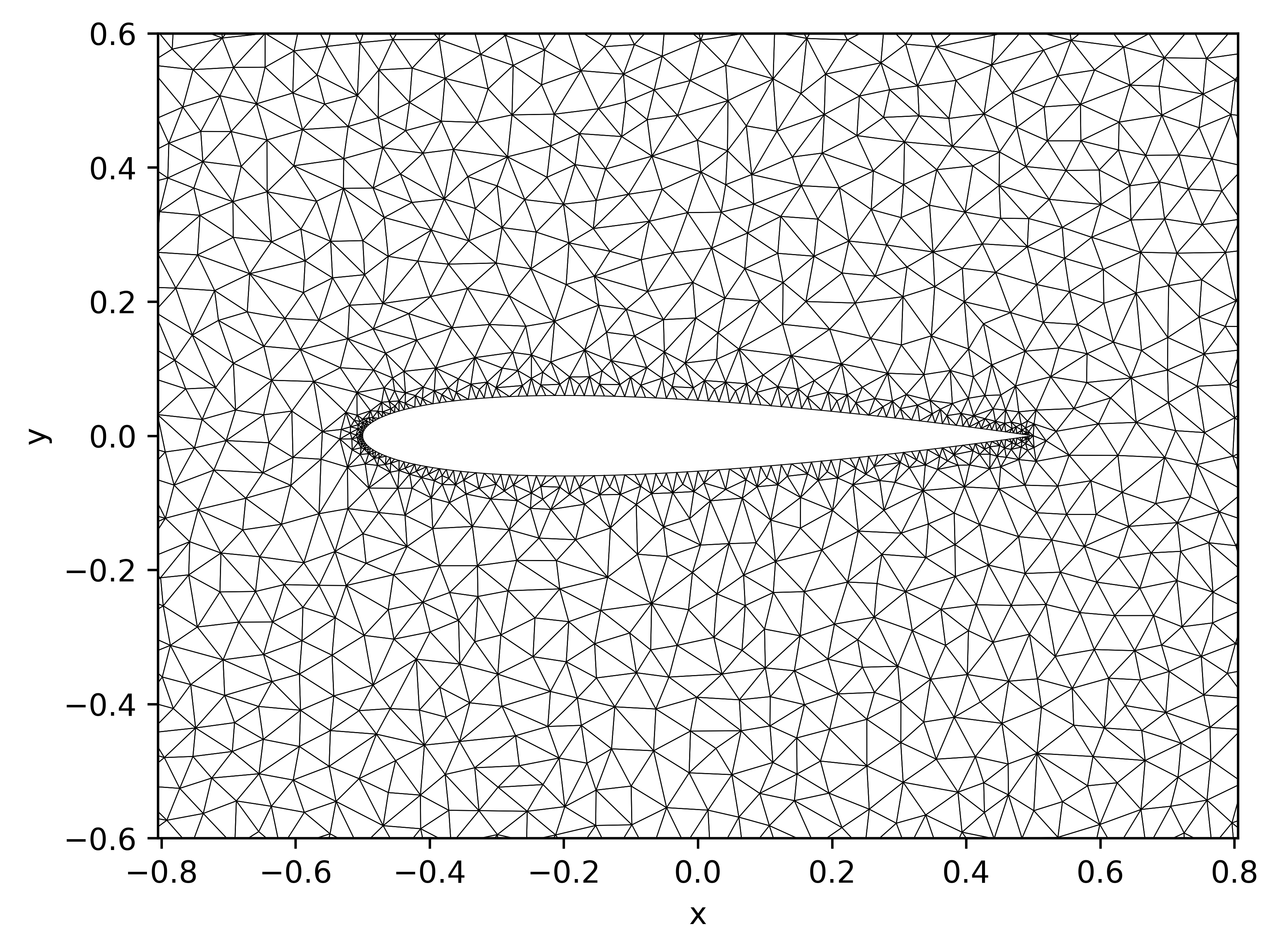}
		\caption{Closeup of the used grid. No refinement next to the airfoil was enforced apart from a minimal inner angle of $28$ degrees per triangle.}
	\end{subfigure}
	\caption{Grid for the NACA 0012 airfoil.}
	\label{fig:Naca0012Grid}
\end{figure}

\begin{figure}
	\begin{subfigure}{\columnwidth}
		\includegraphics[width=\nacawidth]{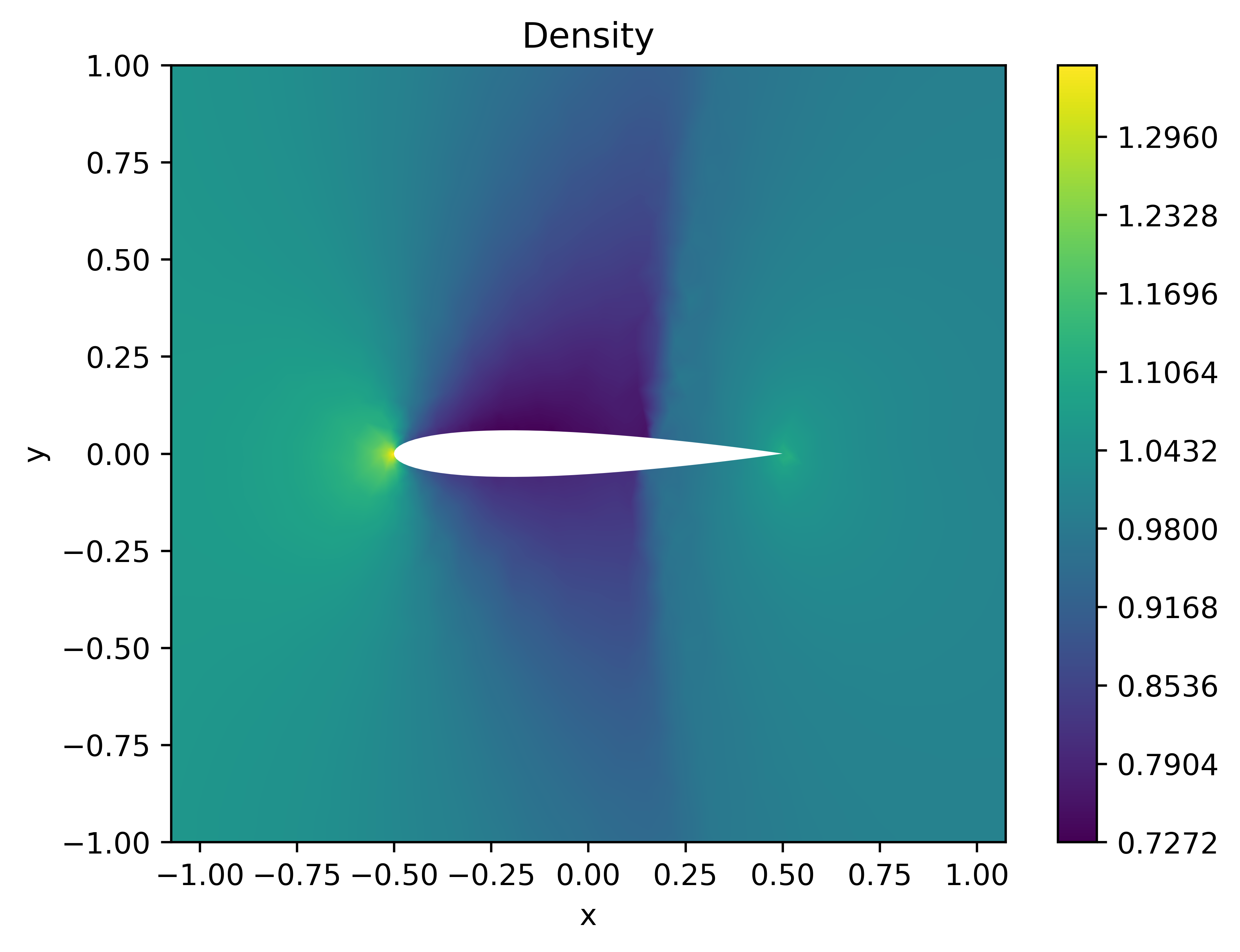}
		\caption{Result for DG1 scheme}
	\end{subfigure}
	\begin{subfigure}{\columnwidth}
		\includegraphics[width=\nacawidth]{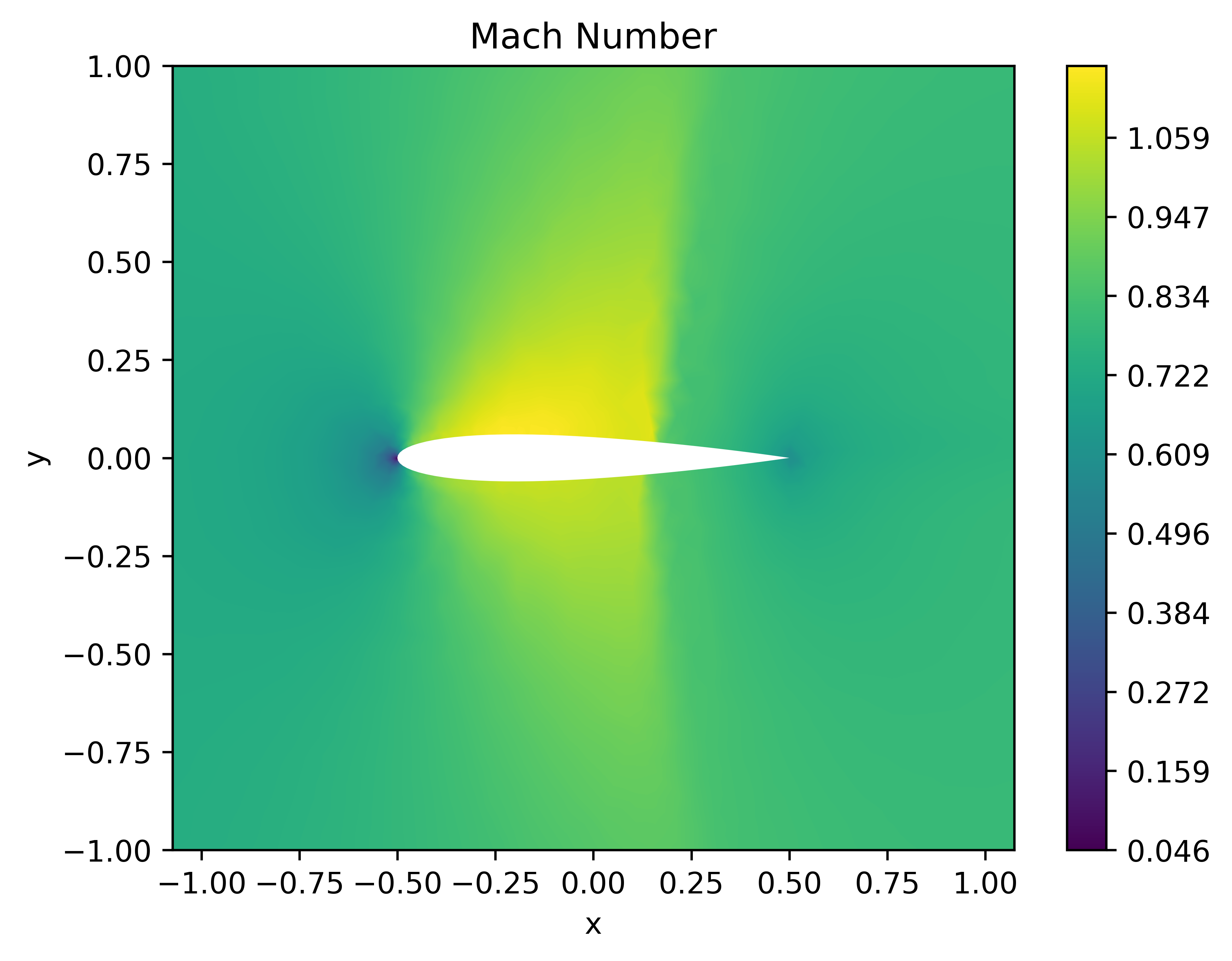}
		\caption{Result for DG1 scheme}
	\end{subfigure}
	\caption{Grid and flow around the NACA 0012 airfoil, $\alpha = 1.25^\circ$, $M_\infty = 0.8$.}
	\label{fig:Naca0012M08}
	\end{figure}

\begin{figure}
	\begin{subfigure}{\columnwidth}
		\includegraphics[width=\nacawidth]{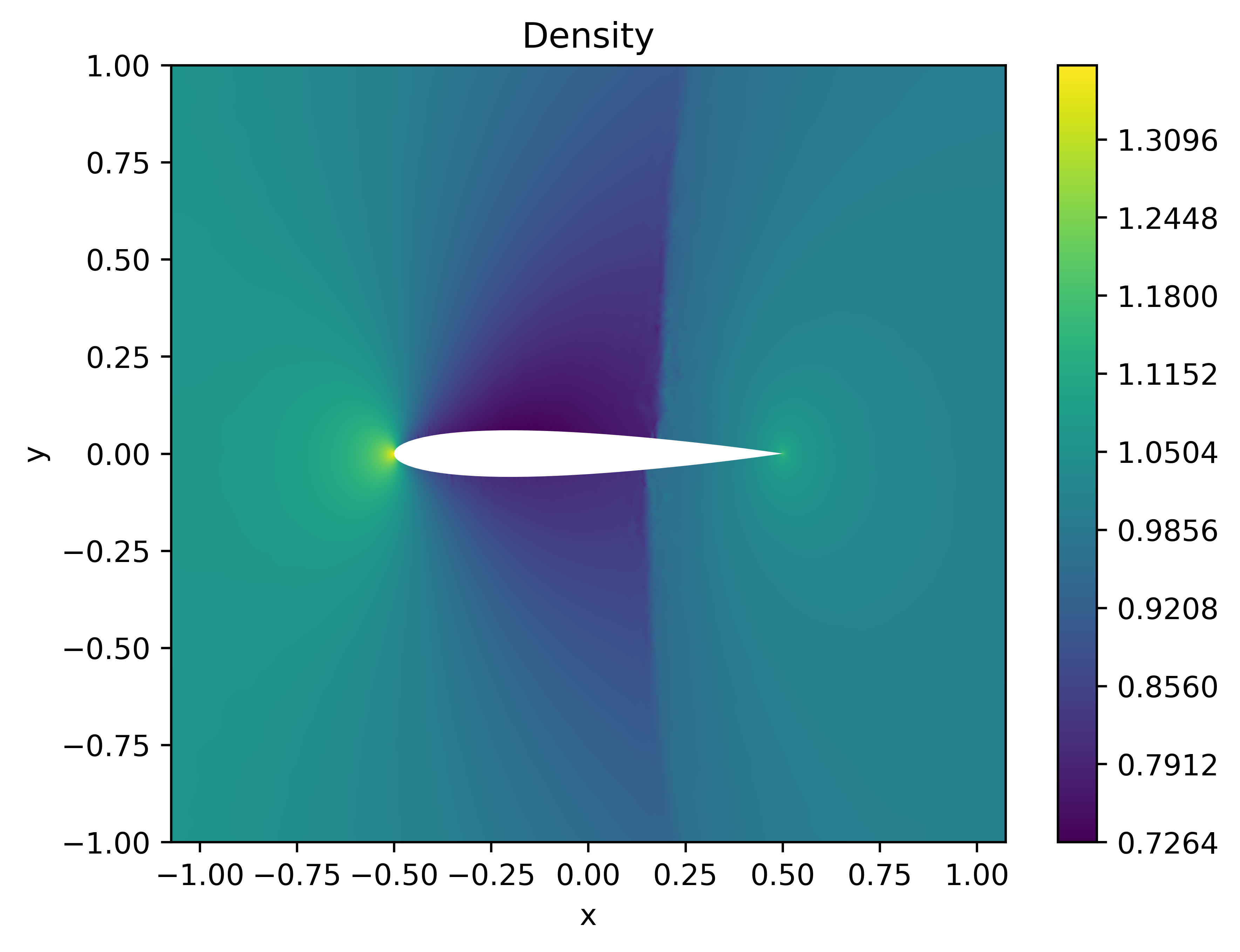}
		\caption{Result for DG3 scheme}
	\end{subfigure}
	\begin{subfigure}{\columnwidth}
		\includegraphics[width=\nacawidth]{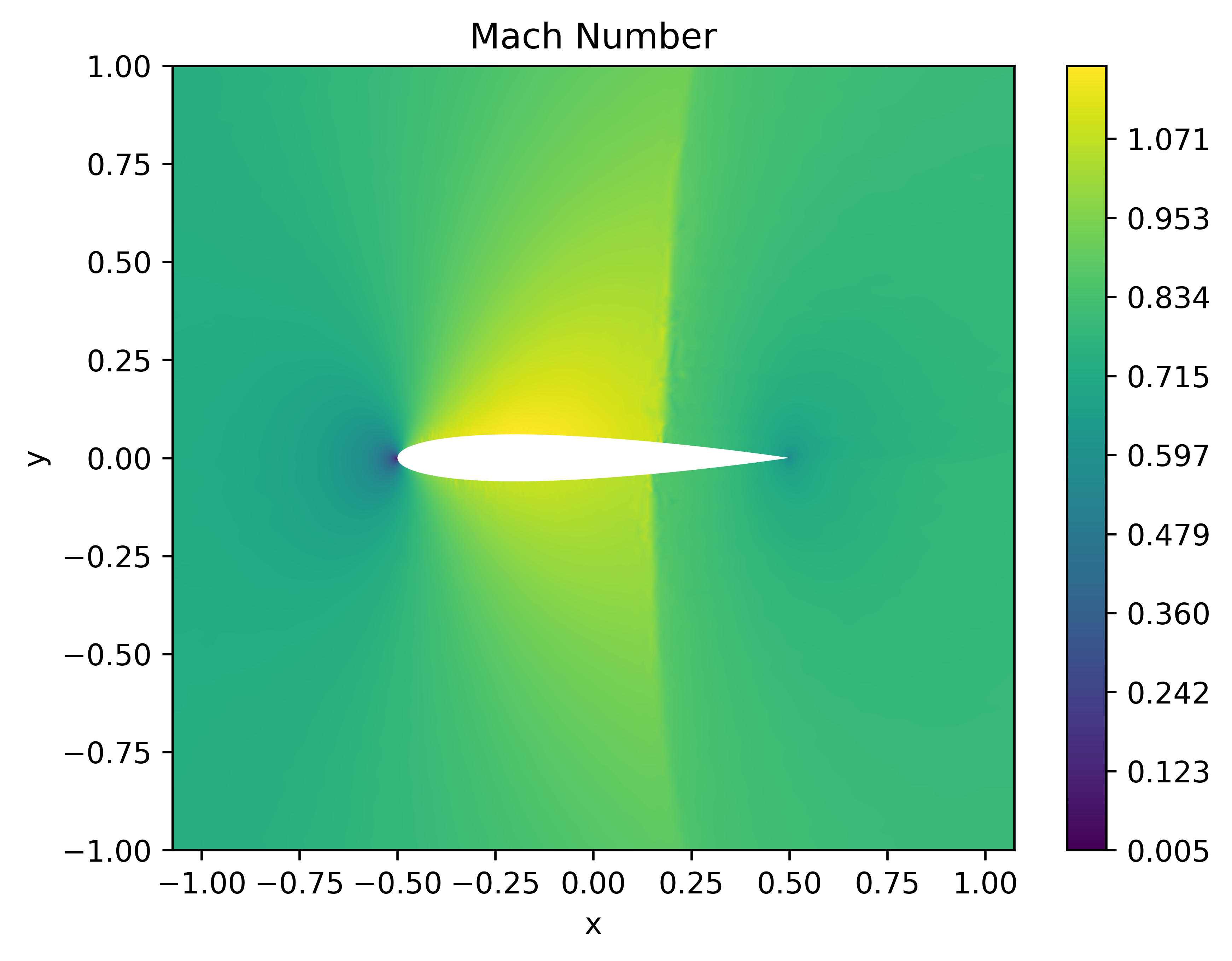}
		\caption{Result for DG3 scheme}
	\end{subfigure}
		\caption{Grid and flow around the NACA 0012 airfoil, $\alpha = 1.25^\circ$, $M_\infty = 0.8$.}
		\label{fig:Naca0012M08p3}
	\end{figure}

\begin{figure}
	\begin{subfigure}{\columnwidth}
		\includegraphics[width=\nacawidth]{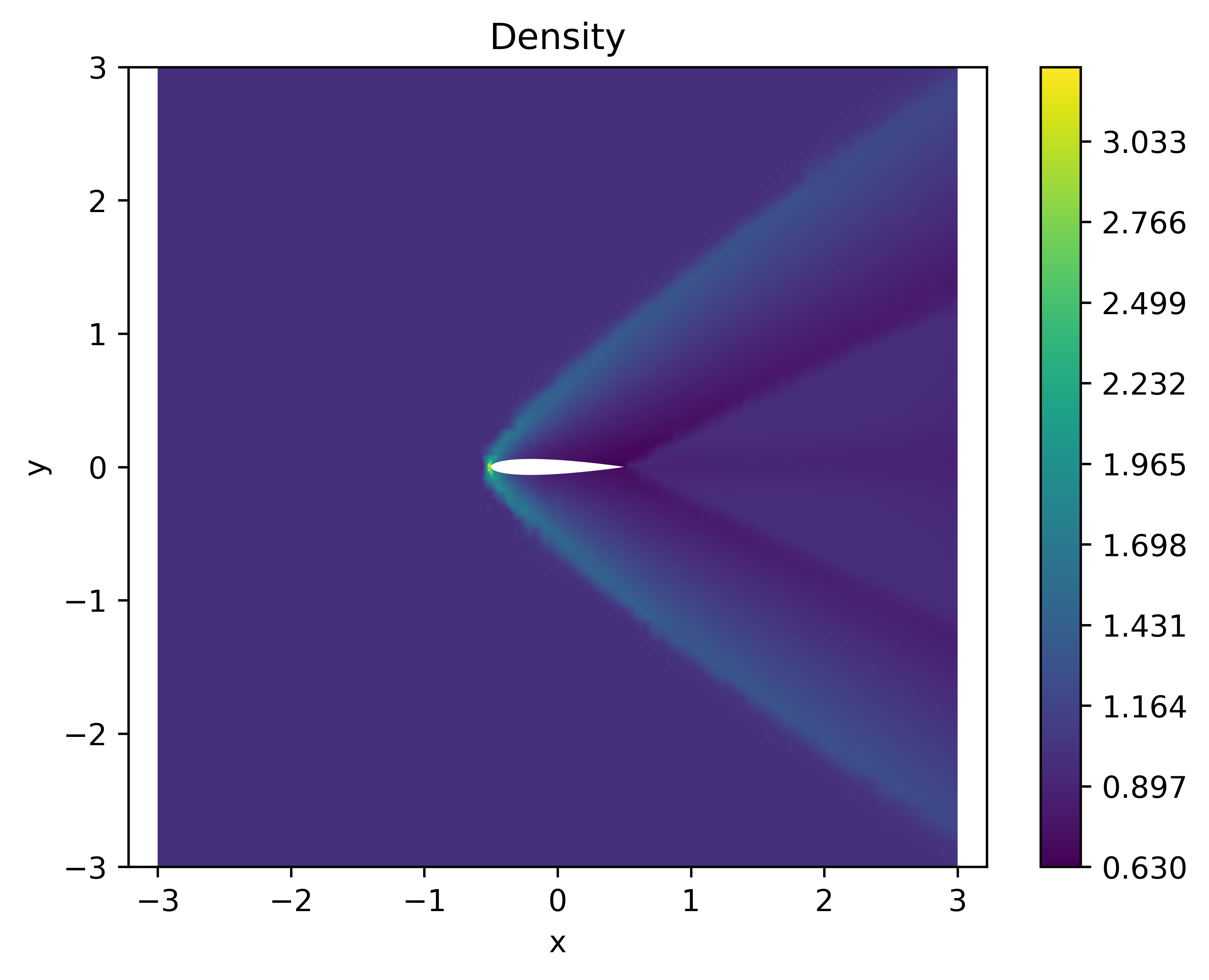}
		\caption{Result for DG1 scheme}
	\end{subfigure}
	\begin{subfigure}{\columnwidth}
			\includegraphics[width=\nacawidth]{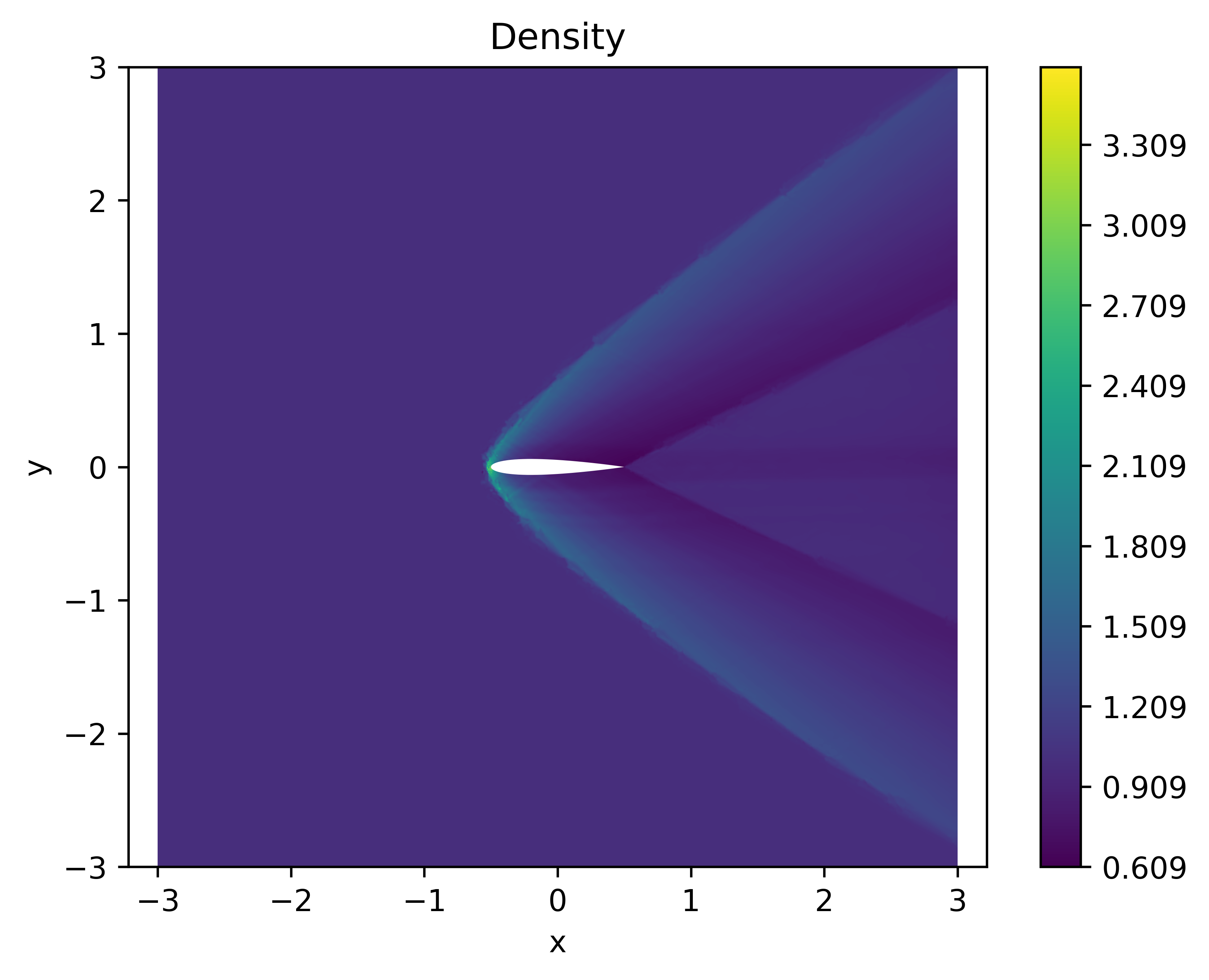}
		\caption{Result for DG3 scheme}
	\end{subfigure}
	\caption{Flow around the NACA 0012 airfoil, $\alpha = 1.25^\circ$, $M_\infty = 2.0$. }
	\label{fig:Naca0012M2}
\end{figure}

\begin{figure}
	\begin{subfigure}{\columnwidth}
		\includegraphics[width=\nacawidth]{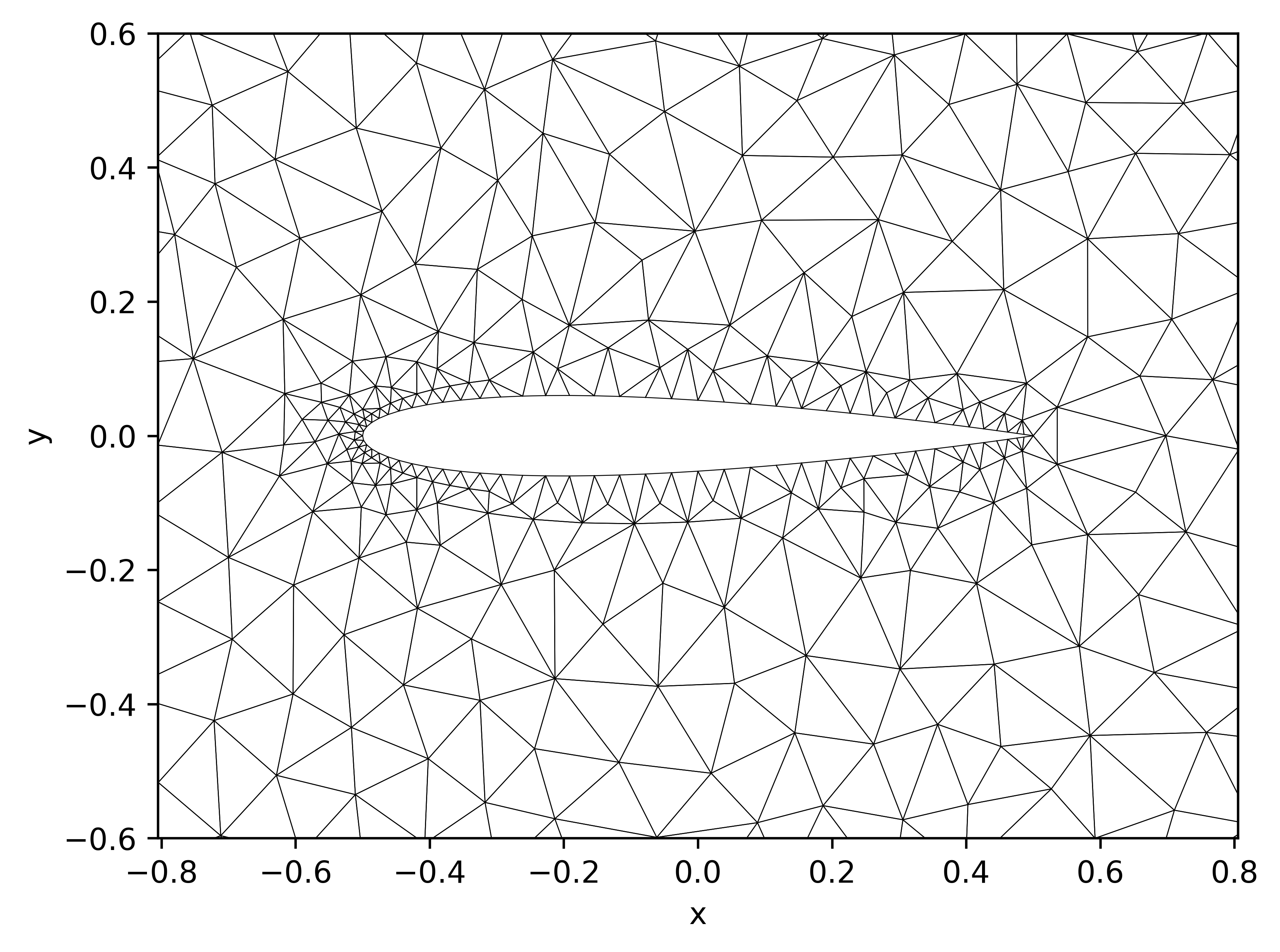}
		\caption{Closeup of the coarse grid of 5844 triangles with a maximal area of $0.01$ per triangle.}
		\end{subfigure}
		\begin{subfigure}{\columnwidth}
			\includegraphics[width=\nacawidth]{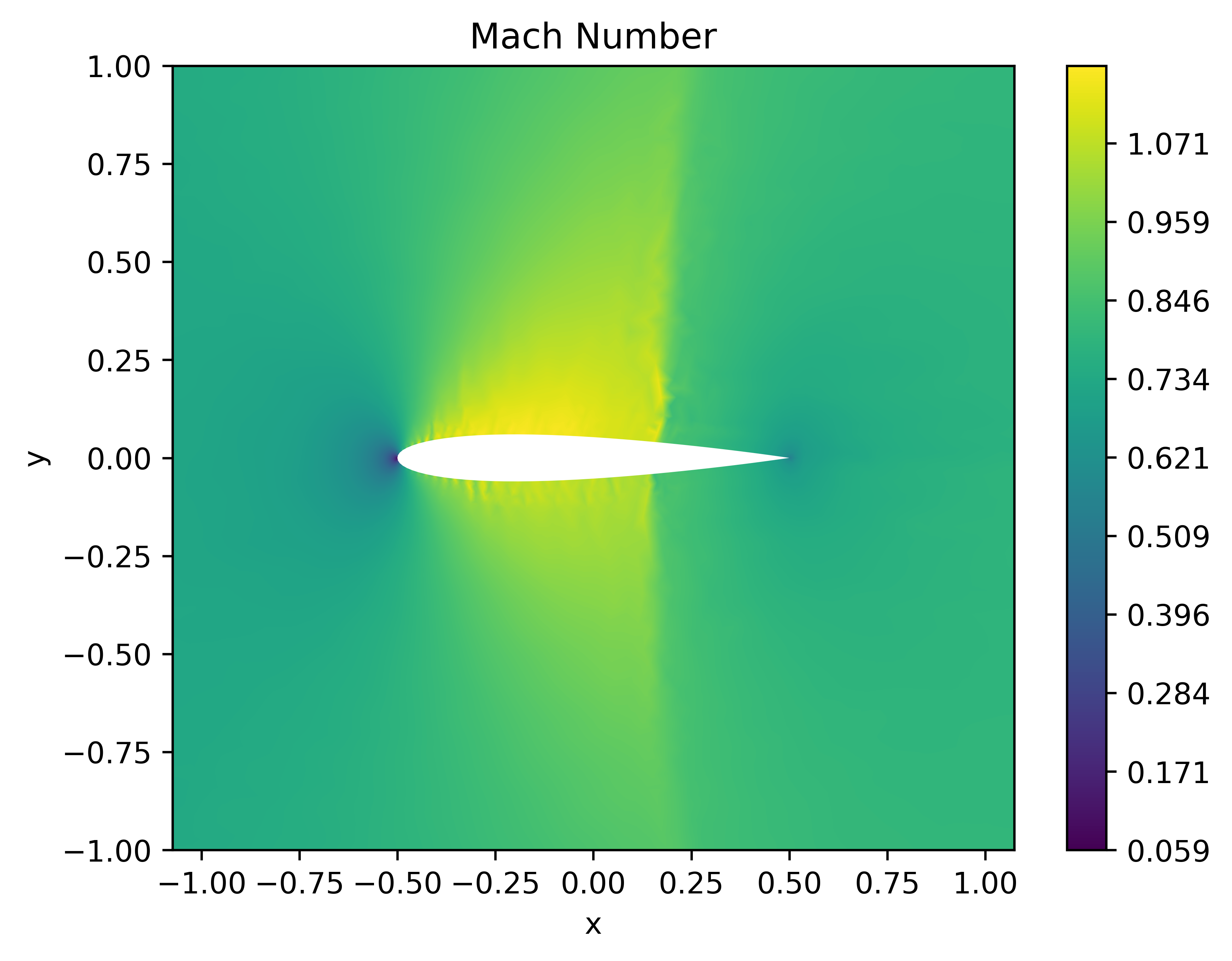}
			\caption{Mach number}
			\end{subfigure}
	\caption{NACA 0012 at low resolution}
	\label{fig:Naca0012M08LR}
	\end{figure}

The NACA 0012 airfoil is a standard benchmark test problem for inviscid flow calculations. The airfoil is tested with an angle of attack of $\alpha = 1.25$ degrees and a free stream velocity of $M_\infty = 0.8$. This results in large areas of smooth flow with two shocks, one on the upper half and one on the lower half of the airfoil \cite{Sonar1997ENO}. This testcase was carried out on a grid with 27771 triangles shown in figure \ref{fig:Naca0012Grid}, and the resulting solution shows the expected shocks in figures \ref{fig:Naca0012M08} and \ref{fig:Naca0012M08p3}. While the high-order calculation shows extremely sharp shocks some oscillations appear at the nose section of the airfoil. We believe that these are not a direct problem of the method but rather a consequence of the piecewise linear boundary representation. Such effects do not appear when the airspeed is increased to $M_\infty = 2$ in figure \ref{fig:Naca0012M2}. The high approximation order can be used to reduce the grid needed to achieve a successful simulation as can be seen in figure \ref{fig:Naca0012M08LR}. There, a low resolution grid with $5844$ triangles was used.
	

	\section{Conclusion} \label{sec:concl}
	The third part of this publication generalizes the entropy rate criterion conforming scheme from \cite{Klein2023StabilizingII} to several space dimensions. 
Apart from an implementation of multidimensional test cases was a detailed justification given to integrate one-dimensional entropy inequality predictors over surfaces of discontinuity to predict the entropy dissipation from multidimensional piecewise smooth functions.
This, together with arguments for the compatibility with first order finite volume schemes, is the basis for the multidimensional implementation of entropy inequality predictors. 
Together with a straight forward generalization of the entropy dissipative conservative filters from the previous part allowed us to implement multidimensional schemes able to calculate solutions of supersonic and transsonic flows essentially oscillation free while achieving high order for smooth flows.
Future developments could include adaptive mesh refinement controlled by the error bound from part one \cite{Klein2023StabilizingI} that is not used in the schemes at the moment, but holds also in several space dimensions.
A generalization of the entropy inequality predictors to finite-volume schemes based on recovery is ongoing work, as to date only error bound based methods are known to the author to enforce the entropy rate criterion in that case \cite{Klein2023ENO}.

	\section{Acknowledgenents} \label{sec:ACK}
	The author would like to thank Thomas Sonar for his comments that improved the manuscript considerably.
	Triangular grids in this work were created using the \guillemetleft Triangle\guillemetright{} program \cite{shewchuk96b}.
		
	\section{Competing Interests}
		The author has no relevant financial or non-financial interests to disclose.
		
	\section{Data Availability}
		The commented implementation of the schemes is available under  \newline \href{https://github.com/simonius/dgdafermos}{https://github.com/simonius/dgdafermos}.
	\section{Bibliography}
	\bibliographystyle{plainnat}
	\bibliography{lit}
\end{document}